%% file: cmp.tex
\documentclass[12pt]{article}
\usepackage{amsmath}
\usepackage{float}
\usepackage{graphicx}
\usepackage{amscd}
\usepackage{amsthm}
\usepackage{amsxtra}
\usepackage{verbatim}
\usepackage[size=2em,heads=vee,Postscript=dvips]{diagrams}
\usepackage{amsfonts}
\usepackage{amssymb}
\usepackage{times}
\diagramstyle[labelstyle=\scriptstyle]
\theoremstyle{plain}
\newtheorem{X}{X}[section]
\newtheorem{theorem}[X]{Theorem}
\newtheorem{proposition}[X]{Proposition}
\newtheorem{lemma}[X]{Lemma}
\newtheorem{corollary}[X]{Corollary}

\newtheorem{E}[X]{}
\theoremstyle{definition}

\newtheorem{definition}[X]{Definition}
\newtheorem{remark}[X]{Remark}
\newtheorem{example}[X]{Example}

\input{defs}

\renewcommand{\Sigma}{\varSigma}
\begin{document}

\title{Abelian Varieties with Complex Multiplication (for Pedestrians)}

\author{J.S. Milne}
\date{}
\maketitle

\setcounter{footnote}{-1}
\begin{abstract}
(June 7, 1998.)\footnote{All footnotes were added in June 1998. TeX style, hence
pagination, changed April 2, 2003.}  This is the text of an article that I wrote and disseminated in 
September 1981, except that I've updated the references, corrected a few 
misprints, and added a table of contents, some footnotes, and an addendum.  

The original article gave a simplified exposition of Deligne's extension of 
the Main Theorem of Complex Multiplication to all automorphisms of the 
complex numbers.  The addendum discusses some additional topics in the 
theory of complex multiplication.  
\end{abstract}

\tableofcontents
\bigskip

The main theorem of Shimura and Taniyama (see Shimura 1971, Theorem 5.15)
describes the action of an automorphism $\tau$ of $\mathbb{C}$ on a polarized
abelian variety of CM-type and its torsion points in the case that $\tau$
fixes the reflex field of the variety. In his Corvallis article (1979),
Langlands made a conjecture concerning conjugates of Shimura varieties that
(see Milne and Shih 1979) leads to a conjectural description of the action of
$\tau$ on a polarized abelian variety of CM-type and its torsion points for
\emph{all} $ \tau$. Recently (July 1981) Deligne proved the
conjecture\footnote{That is, the conjectural description of the action of
$\tau$ on a polarized abelian variety of CM-type, not Langlands's conjecture!}
(Deligne 1982b). Deligne expresses his result as an identity between two
pro-reductive groups, the Taniyama group of Langlands and his own motivic
Galois group associated with the Tannakian category of motives of abelian
varieties over $\mathbb{Q}$ of potential CM-type. Earlier (April 1981) Tate
gave a more down-to-earth conjecture than that stated in (Milne and Shih 1979)
and partially proved his conjecture (Tate 1981).

The purpose of these notes is to use Deligne's ideas to give as direct a proof
as possible of the conjecture in the form stated by Tate. It is also checked
that the three forms of the conjecture, those in Deligne 1982b, Milne and Shih
1979, and Tate 1981 are compatible. Also, Tate's ideas are used to simplify
the construction of the Taniyama group. In the first three sections, I have
followed Tate's manuscript (1981) very closely, sometimes word-for-word.

These notes are a rough write-up of two of my lectures at the conference on
Shimura Varieties, Vancouver, 17--25 August, 1981. In the remaining lectures I
described how the result on abelian varieties of CM-type could be applied to
give a proof of Langlands's conjecture on conjugates of Shimura varieties for
most\footnote{In fact all, see Milne 1983.} (perhaps all) Shimura varieties.

\subsubsection{Notations}

We let $\hat{\mathbb{Z}}=\plim{\mathbb{Z}} /m\mathbb{Z}$ and $\mathbb{A}%
_{f}=\hat{\mathbb{Z}}\otimes\mathbb{Q}$. For a number field $E$,
$\mathbb{A}_{f,E}=\mathbb{A}_{f}\otimes E$ is the ring of finite ad\`eles and
$\mathbb{A}_{E}=\mathbb{A}_{f,E}\times(E\otimes\mathbb{R})$ the full ring of
ad\`eles. When $ E$ is a subfield of $\mathbb{C}$, $E^{\ab}$ and $E^{\al}$
denote respectively the maximal abelian extension of $ E$ in $\mathbb{C}$ and
the algebraic closure of $E$ in $\mathbb{C}$. Complex conjugation is denoted
by $ \iota$.

For a number field $E$, $\rec_{E}:\mathbb{A}_{E}^{\times}\to\Gal(E^{\ab} /E)$
is the reciprocity law, normalized so that a prime element parameter
corresponds to the inverse of the usual (arithmetic) Frobenius: if
$a\in\mathbb{A}_{f,E}^{\times}$ has $ v$-component a prime element $a_{v}$ in
$E_{v}$ and $w$-component $a_{w}=1$ for $w\neq v$, then $\rec_{E}%
(a)=\sigma^{-1}$ if $\sigma x\equiv x^{\mathbb{N}(v)}\mod {\mathfrak{p}}_{v}$. When
$E$ is totally complex, $\rec_{E}$ factors into $\mathbb{A}_{E}^{\times}%
\to\mathbb{A}_{f,E}^{\times}\xr {r_E}\Gal(E^{a b}/E).$ The cyclotomic
character $\chi=\chi_{\text{\textrm{cyc}}}:\Aut (\mathbb{C})\to\hat
{\mathbb{Z}}^{\times}\subset\mathbb{A}_{f}^{\times}$ is the homomorphism such
that $\tau\zeta=\zeta^{\chi(\tau)}$ for every root of $ 1$ in $\mathbb{C}$.
The composite $r_{E}\circ\chi=\text{\textrm{Ver}}_{E/\mathbb{Q}}$, the
Verlagerung map $\Gal(\mathbb{Q}^{\al}/\mathbb{Q})^{\ab}\to\Gal(\mathbb{Q}%
^{\al}/E)^{\ab}$.

When $T$ is a torus over $E$, $X_{*}(T)$ is the cocharacter group $\Hom_{
E^{\al}}(\mathbb{G}_{m},T)$ of $T$.

*Be wary\footnote{This is universally good advice, but I believe the signs
here to be correct.} of signs.*

\section{Statement of the Theorem}

Let $A$ be an abelian variety over $ \mathbb{C}$, and let $K$ be a subfield of
$\End(A)\otimes\mathbb{Q}$ of degree $2\dim A$ over $\mathbb{Q}$. The
representation of $K$ on the tangent space to $A$ at zero is of the form $
\oplus_{\phi\in\Phi}\phi$ with $\Phi$ a subset of $\Hom(K,\mathbb{C})$. A
\emph{Riemann form\/} for $ A$ is a $\mathbb{Q}$-bilinear skew-symmetric form
$\psi$ on $H_{1}(A,\mathbb{Q})$ such that
\[
(x,y)\mapsto\psi(x,iy):H_{1}(A,\mathbb{R})\times H_{1}(A,\mathbb{R}%
)\to\mathbb{R}%
\]
is symmetric and positive definite. We assume that there exists a Riemann form
$\psi$ compatible with the action of $K$ in the sense that
\[
\psi(ax,y)=\psi(x,(\iota a)y),\quad a\in K,\quad x,y\in H_{1}(A ,\mathbb{Q}).
\]
Then $K$ is a CM-field, and $\Phi$ is a CM-type on $K$, i.e., $\Hom
(K,\mathbb{C})=\Phi\cup\iota\Phi$ (disjoint union). The pair
$(A,K\hookrightarrow\End(A)\otimes\mathbb{Q} )$ is said to be of
\emph{CM-type} $(K,\Phi)$. For simplicity, we assume that $K\cap
\End(A)=\mathcal{O}_{K}$, the full ring of integers in $K$.

Let $\mathbb{C}^{\Phi}$ be the set of complex-valued functions on $ \Phi
$\emph{,\/} and embed $K$ into $\mathbb{C}^{\Phi}$ through the natural map
$a\mapsto(\phi(a))_{\phi\in\Phi}$. There then exist a $ \mathbb{Z}$-lattice
$\mathfrak{a}$ in $K$ stable under $\mathcal{O}_{K}$, an element $t\in
K^{\times}$, and an $ \mathcal{O}_{K}$-linear analytic isomorphism
$\theta:\mathbb{C}^{\Phi}/\mathfrak{a}\to A$ such that $\psi(x,y)=\Tr_{
K/\mathbb{Q}}(tx\cdot\iota y)$ where, in the last equation, we have used
$\theta$ to identify $H_{1}(A,\mathbb{Q})$ with $\mathfrak{a}\otimes
\mathbb{Q}=K$. The variety is said to be of \emph{type} $(K,\Phi
;\mathfrak{a},t)$ relative\footnote{See Shimura 1971, pp 126--129.} to $
\theta$. The type determines the triple $(A,K\hookrightarrow\End(A)\otimes
\mathbb{Q} ,\psi)$ up to isomorphism. Conversely, the triple determines the
type up to a change of the following form: if $\theta$ is replaced by
$\theta\circ a^{-1}$, $a\in K^{ \times}$, then the type becomes $(K,\Phi
;a\mathfrak{a},\frac t{a\cdot\iota a})$.

Let $\tau\in\Aut(\mathbb{C})$. Then $K\hookrightarrow\End(A)\otimes\mathbb{Q}$
induces a map $K\hookrightarrow\End(\tau A)\otimes\mathbb{Q}$, so that $\tau
A$ also has complex multiplication by $K$. The form $\psi$ is associated with
a divisor $D$ on $A$, and we let $\tau\psi$ be the Riemann form for $ \tau A$
associated with $\tau D$. It has the following characterization: after
multiplying $ \psi$ with a nonzero rational number, we can assume that it
takes integral values on $H_{1}(A,\mathbb{Z})$; define $\psi_{m}$ to be the
pairing $A_{m}\times A_{m}\to\mu_{m}$, $(x,y)\mapsto\exp(\frac{2\pi i\cdot
\psi(x,y)}m)$; then $(\tau\psi)_{m}(\tau x,\tau y)=\tau(\psi_{m}(x,y))$.

In the next section we shall define (following Tate) for each CM-type $
(K,\Phi)$ a map $f_{\Phi}:\Aut(\mathbb{C})\to\mathbb{A}_{f,K}^{\times
}/K^{\times}$ such that
\[
f_{\Phi}(\tau)\cdot\iota f_{\Phi}(\tau)=\chi(\tau)K^{\times} ,\quad
\text{\textrm{all }}\tau\in\Aut(\mathbb{C}).
\]
We can now state the new main theorem of complex multiplication in the form
first appearing (as a conjecture) in Tate 1981.

\begin{theorem}[Shimura, Taniyama, Langlands, Deligne]
Suppose $A$ has type $(K,\Phi ;{\mathfrak{a}},t)$ relative to $\theta
:\mathbb{C}^{\Phi}/{\mathfrak{a}}\xr {\approx}A$.  Let $\tau\in\Aut(\mathbb{C}
)$, and
let $f\in \mathbb{A}_{f,K}^{\times}$ lie in $f_{\Phi}(\tau )$.
\begin{enumerate}
\item The variety $\tau A$ has type
\[(K,\tau\Phi ;f{\mathfrak{a}},\frac {t\chi (\tau )}{f\cdot\iota f})\]
relative to $\theta'$ say.
\item It is possible to choose $\theta'$ so that
\[\begin{diagram}
{\mathbb{A}}_{f,K}&\rOnto&{\mathbb{A}}_{f,K}
/{\mathfrak{a}}\otimes\hat{{\mathbb{Z}}}\cong{K}/{\mathfrak{a}}&\rTo{\theta}&A_{\text{tors}}\\
\dTo{}{f}&&&&\dTo{}{\tau}\\
{\mathbb{A}}_{f,K}&\rOnto&{\mathbb{A}}_{f,K}
/(f{\mathfrak{a}}\otimes\hat{{\mathbb{Z}}})\cong{K}/f{\mathfrak{a}}&\rTo{\theta'}&\tau{}A_{\text{tors}}
\end{diagram}\]
commutes, where $A_{\text{\rm tors}}$ denotes the torsion subgroup of $
A$.
\end{enumerate}
\end{theorem}

\begin{remark}
Prior to its complete proof, the theorem was known in
three\footnote{Shimura (1977) investigated the question in some further
special cases.  After explaining that the action of a general automorphism
of $\mathbb{C}$ on an elliptic curve of CM-type can be obtained from knowing the
actions of complex conjugation and those automorphisms fixing its reflex
field, he concludes rather pessimistically that ``In the higher-dimensional
case, however, no such general answer seems possible.''} important cases.
\begin{enumerate}
\item If $\tau$ fixes the reflex field of $(K,\Phi )$, then the theorem becomes the old
main theorem of complex multiplication, proved by Shimura and Taniyama
(see (2.7) below).  This case is used in the proof of the general result.
\item\ Tate (1981) proved part (a) of the theorem, and he showed that (b)
holds when $f$ is replaced by $fe$, some $e\in \mathbb{A}_{f,K_0}^{
\times}$ with $e^2=1$, where $K_0$ the
maximal real subfield of $K$.  We include Tate's proof of his result, although
it is not necessary for the general case.
\item\ Shih (1976) proved the theorem under the assumption that there exists
an automorphism $\sigma$ of $K$ of order $2$ such that $\tau (\Phi
\cap\Phi\sigma )=\Phi\cap\Phi\sigma$ and
$\tau (\Phi\cap\iota\Phi\sigma )=\Phi\cap\iota\Phi\sigma$ for all automorphisms $
\tau$ of $\mathbb{C}$.  As we shall see, his
proof is a special case of the general proof.
\end{enumerate}
\end{remark}

We now restate the theorem in more invariant form. Let
\[
TA\overset{\text{\textrm{df}}}{=}\plim A_{m}(\mathbb{C})\cong\plim(\frac
1mH_{1}(A,\mathbb{Z})/H_{1}(A,\mathbb{Z}))=H_{1}(A,\hat{\mathbb{Z}})
\]
(limit over all positive integers $m$), and let
\[
V_{f}A=TA\otimes_{\mathbb{Z}}\mathbb{Q}=H_{1}(A,\mathbb{Q})\otimes
_{\mathbb{Q}} \mathbb{A}_{f}.
\]
Then $\psi$ gives rise to a pairing
\[
\psi_{f}=\plim\psi_{m}:V_{f}A\times V_{f}A\to\mathbb{A}_{f}(1)
\]
where $\mathbb{A}_{f}(1)=(\plim\mu_{m}(\mathbb{C}))\otimes\mathbb{Q}$.

\begin{theorem}
Let $A$ have type $(K,\Phi )$; let $\tau\in\Aut(\mathbb{C})$, and let $
f\in f_{\Phi}(\tau )$.
\begin{enumerate}
\item$\tau A$ is of type $(K,\tau\Phi )$;
\item there is an $K$-linear isomorphism $\alpha :H_1(A,\mathbb{Q})
\to H_1(\tau A,\mathbb{Q})$ such that
\begin{enumerate}
\item$\psi (\frac {\chi (\tau )}{f\cdot\iota f}x,y)=(\tau\psi )(\alpha
x,\alpha y),\quad x,y\in H_1(A,\mathbb{Q})$;
\item\ the\footnote{Note that both $f\in \mathbb{A}_{f,K}^{\times}$ and the $
K$-linear isomorphism $\alpha$
are uniquely determined up to multiplication by an element of $K^{
\times}$.  Changing
the choice of one changes that of the other by the same factor.}  diagram
\[\begin{diagram}[heads=LaTeX]
V_f(A)&\rTo{f}&V_f(A)\\
&\rdTo{\tau}&\dTo{}{\alpha\otimes1}\\
&&V_f(\tau{A})\end{diagram}\]
commutes.
\end{enumerate}
\end{enumerate}
\end{theorem}

\begin{lemma}
The statements (1.1) and (1.3) are equivalent.
\end{lemma}
\begin{proof}
Let $\theta$ and $\theta'$ be as in (1.1), and let $\theta_1:K\xr {\approx}
H_1(A,\mathbb{Q})$ and $\theta_1':K\xr {\approx}H_1(\tau A,\mathbb{Q}
)$ be the
$K$-linear isomorphisms induced by $\theta$ and $\theta'$.  Let $
\chi =\chi (\tau )/f\cdot\iota f$ --- it is an
element of $K^{\times}$.  Then
\begin{eqnarray*}
\psi (\theta_1(x),\theta_1(y))&=&\Tr(tx\cdot\iota y)\\
(\tau\psi )(\theta_1'(x),\theta_1'(y))&=&\Tr(t\chi x\cdot\iota y)\end{eqnarray*}
and
\[\begin{CD}
{\mathbb{A}}_{f,K}@>{\theta_1}>>V_f(A)\\
@VV{f}V@VV{\tau}V\\
{\mathbb{A}}_{f,K}@>\theta_1'>>V_f(\tau{A})
\end{CD}\]
commutes.  Let $\alpha =\theta_1'\circ\theta_1^{-1}$; then
\[\tau\psi (\alpha x,\alpha y)=\Tr(t\chi\theta_1^{-1}(x)\cdot\iota
\theta_1^{-1}(y))=\psi (\chi x,y)\]
and (on $V_f(A)$),
\[\tau =\theta_1'\circ f\circ\theta_1^{-1}=\theta_1'\circ\theta_1^{
-1}\circ f=\alpha\circ f.\]
Conversely, let $\alpha$ be as in (1.3) and choose $\theta_1'$ so that $
\alpha =\theta_1'\circ\theta_1^{-1}$.  It is then
easy to check (1.1).
\end{proof}

\section{Definition of $f_{\Phi}(\tau)$}

Let $(K,\Phi)$ be a CM-type. Choose an embedding $K\hookrightarrow\mathbb{C},$
and extend it to an embedding $i:K^{\ab}\hookrightarrow\mathbb{C}$. Choose
elements $ w_{\rho}\in\Aut(\mathbb{C})$, one for each $\rho\in
\Hom(K,\mathbb{C})$, such that
\[
w_{\rho}\circ i|K=\rho,\quad w_{\iota\rho}=\iota w_{\rho}.
\]
For example, choose $w_{\rho}$ for $\rho\in\Phi$ (or any other CM-type) to
satisfy the first equation, and then define $w_{\rho}$ for the remaining
$\rho$ by the second equation. For any $\tau\in\Aut(\mathbb{C})$, $w_{\tau
\rho}^{-1}\tau w_{\rho}\circ i|K=w_{\tau\rho}^{-1}\circ\tau\rho|K=i$. Thus
$i^{-1}\circ w_{\tau\rho}^{-1}\tau w_{\rho}\circ i\in\Gal(K^{\ab} /K)$, and we
can define $F_{\Phi}:\Aut(\mathbb{C})\to\Gal(K^{\ab}/K )$ by
\[
F_{\Phi}(\tau)=\prod_{\phi\in\Phi}i^{-1}\circ w_{\tau\phi}^{-1} \tau w_{\phi
}\circ i.
\]

\begin{lemma}
The element $F_{\Phi}$ is independent of the choice of $\{w_{\rho}
\}$.
\end{lemma}
\begin{proof}
Any other choice is of the form $w_{\rho}'=w_{\rho}h_{\rho}$, $h_{
\rho}\in\Aut(\mathbb{C}/iK)$.  Thus $F_{\Phi}(\tau )$ is
changed by $i^{-1}\circ (\prod_{\phi\in\Phi}h_{\tau\phi}^{-1}h_{\phi}
)\circ i$.  The conditions on $w$ and $w'$ imply that
$h_{\iota\rho}=h_{\rho}$, and it follows that the inside product is $
1$ because $\tau$ permutes the
unordered pairs $\{\phi ,\iota\phi \}$ and so $\prod_{\phi\in\Phi}
h_{\phi}=\prod_{\phi\in\Phi}h_{\tau\phi}$.
\end{proof}

\begin{lemma}
The element $F_{\Phi}$ is independent of the choice of $i$ (and $
K\hookrightarrow \mathbb{C}$).
\end{lemma}
\begin{proof}
Any other choice is of the form $i'=\sigma\circ i$, $\sigma\in\Aut
(\mathbb{C})$.  Take $w_{\rho}'=w_{\rho}\circ\sigma^{-1}$,
and then
\[F_{\Phi}'(\tau )=\prod i^{\prime -1}\circ (\sigma w_{\tau\phi}^{
-1}\tau w_{\phi}\sigma^{-1})\circ i'=F_{\Phi}(\tau ).\]
\end{proof}

Thus we can suppose $K\subset\mathbb{C}$ and ignore $i$; then
\[
F_{\Phi}(\tau)=\prod_{\phi\in\Phi}w_{\tau\phi}^{-1}\tau w_{\phi}\mod\Aut
(\mathbb{C}/K^{ab})
\]
where the $w_{\rho}$ are elements of $\Aut(\mathbb{C})$ such that
\[
w_{\rho}|K=\rho,\quad w_{\iota\rho}=iw_{\rho}.
\]

\begin{proposition}
For any $\tau\in\Aut(\mathbb{C})$, there is a unique $f_{\Phi}(\tau
)\in \mathbb{A}_{f,K}^{\times}/K^{\times}$ such that
\begin{enumerate}
\item\ $r_K(f_{\Phi}(\tau ))=F_{\Phi}(\tau )$;
\item\ $f_{\Phi}(\tau )\cdot\iota f_{\Phi}(\tau )=\chi (\tau )K^{
\times}$, $\chi =\chi_{\text{\rm cyc}}$.
\end{enumerate}
\end{proposition}
\begin{proof}
Since $r_K$ is surjective, there is an $f\in \mathbb{A}_{f,K}^{\times}
/K^{\times}$ such that $r_K(f)=F_{\Phi}(\tau )$.
We have
\begin{eqnarray*}
r_K(f\cdot\iota f)&=&r_K(f)\cdot r_K(\iota f)\\
&=&r_K(f)\cdot\iota r_K(f)\iota^{-1}\\
&=&F_{\Phi}(\tau )\cdot F_{\iota\Phi}(\tau )\\
&=&V_{K/\mathbb{Q}}(\tau ),\end{eqnarray*}
where $V_{K/\mathbb{Q}}:\Gal(\mathbb{Q}^{\al}/\mathbb{Q})^{\text{\rm ab}}
\to\Gal(\mathbb{Q}^{\al}/K)^{\text{\rm ab}}$ is the transfer (Verlagerung) map.
As $V_{K/\mathbb{Q}}=r_K\circ\chi$, it follows that $f\cdot\iota f=
\chi (\tau )K^{\times}\mod(\Ker r_K)$.  The next
lemma shows that $1+\iota$ acts bijectively on $\Ker(r_K)$, and so there is a
unique $a\in\Ker r_K$ such that $a\cdot\iota a=(f\cdot\iota f/\chi
(\tau ))K^{\times}$; we must take
$f_{\Phi}(\tau )=f/a$.
\end{proof}

\begin{lemma}
The kernel of $r_K:\mathbb{A}_{f,K}^{\times}/K^{\times}\to\Gal(K^{\text {\rm ab}}
/K)$ is uniquely divisible by all integers,
and its elements are fixed by $\iota$.
\end{lemma}
\begin{proof}
The kernel of $r_K$ is $\overline {K^{\times}}/K^{\times}$ where $\overline {
K^{\times}}$ is the closure of $K^{\times}$ in $\mathbb{A}_{f,K}^{\times}$.  It is
also equal to $\bar {U}/U$ for any subgroup $U$ of ${\cal O}_K^{\times}$ of finite index.  A theorem of
Chevalley (see Serre 1964, 3.5) shows that $\mathbb{A}_{f,K}^{\times}$ induces the pro-finite
topology on $U$.  If we take $U$ to be contained in the real subfield of $
K$ and
torsion-free, then it is clear that $\bar {U}/U$ is fixed by $\iota$ and (being isomorphic
to $(\hat {\mathbb{Z}}/\mathbb{Z})^{\dim A}$) uniquely divisible.
\end{proof}

\begin{remark}
A more direct definition of $f_{\Phi}(\tau )$, but one involving the Weil group, can be
found in (7.2).
\end{remark}

\begin{proposition}
The maps $f_{\Phi}:\Aut(\mathbb{C})\to\mathbb{A}_{f,K}^{\times}/K^{\times}$ have the following properties:
\begin{enumerate}
\item\ $f_{\Phi}(\sigma\tau )=f_{\tau\Phi}(\sigma )\cdot f_{\Phi}
(\tau )$;
\item\ $f_{\Phi (\tau^{-1}|K)}(\sigma )=\tau f_{\Phi}(\sigma )$ if $
\tau K=K$;
\item\ $f_{\Phi}(\iota )=1$.
\end{enumerate}
\end{proposition}
\begin{proof}
Let $f=f_{\tau\Phi}(\sigma )\cdot f_{\Phi}(\tau )$.  Then
\[r_K(f)=F_{\tau\Phi}(\sigma )\cdot F_{\Phi}(\tau )=\prod_{\phi\in
\Phi}w_{\sigma\tau\phi}^{-1}\sigma w_{\tau\phi}w_{\tau\phi}^{-1}\tau
w_{\phi}=F_{\Phi}(\sigma\tau )\]
and $f\cdot\iota f=\chi (\sigma )\chi (\tau )K^{\times}=\chi (\sigma
\tau )K^{\times}$.  Thus $f$ satisfies the conditions that
determine $f_{\Phi}(\sigma\tau )$.  This proves (a), and (b) and (c) can be proved similarly.
\end{proof}

Let $E$ be the reflex field for $(K,\Phi)$, so that $\Aut(\mathbb{C}
/E)=\{\tau\in\Aut(\mathbb{C})\mid\tau\Phi=\Phi\}$. Then $\Phi\Aut(\mathbb{C}%
/K)\overset{\text{\textrm{df}}}{=}\cup_{\phi\in\Phi}\phi\cdot\Aut(\mathbb{C}%
/K)$ is stable under the left action of $\Aut (\mathbb{C}/E)$, and we write
\[
\Aut(\mathbb{C}/K)\Phi^{-1}=\cup\psi\Aut(\mathbb{C}/E)\qquad
(\text{\textrm{disjoint union}} );
\]
the set $\Psi=\{\psi|E\}$ is a CM-type for $E$, and $(E,\Psi)$ is the reflex
of $ (K,\Phi)$. The map $a\mapsto\prod_{\psi\in\Psi}\psi(a):E\to\mathbb{C}$
factors through $ K$ and defines a morphism of algebraic tori $\Psi^{\times
}:E^{\times}\to K^{\times}$. The (old) main theorem of complex multiplication
states the following: let $\tau\in\Aut(\mathbb{C}/E)$, and let $a\in
\mathbb{A}_{f,E}^{\times}/E^{\times}$ be such that $r_{E}(a)=\tau$; then (1.1)
is true after $f$ has been replaced by $ \Psi^{\times}(a)$. (See Shimura 1971,
Theorem 5.15; the sign differences result from different conventions for the
reciprocity law and the actions of Galois groups.) The next result shows that
this is in agreement with (1.1).

\begin{proposition}
For any $\tau\in\Aut(\mathbb{C}/E)$ and $a\in\mathbb{A}_{f,E}^{\times}
/E^{\times}$ such that $r_E(a)=\tau$, $\Psi^{\times}(a)\in f_{\Phi}
(\tau )$.
\end{proposition}
\begin{proof}
Partition $\Phi$ into orbits, $\Phi =\cup_j\Phi_j$, for the left action of $\Aut
(\mathbb{C}/E)$.  Then
$\Aut(\mathbb{C}/K)\Phi^{-1}=\cup_j\Aut(\mathbb{C}/K)\Phi_j^{-1}$, and
\[\Aut(\mathbb{C}/K)\Phi_j^{-1}=\Aut(\mathbb{C}/K)(\sigma_j^{-1}\Aut(
\mathbb{C}/E))=(\Hom_K(L_j,\mathbb{C})\circ\sigma_j^{-1})\Aut(\mathbb{C}
/E)\]
where $\sigma_j$ is any element of $\Aut(\mathbb{C})$ such that $\sigma_
j|K\in\Phi_j$ and $L_j=(\sigma_j^{-1}E)K$.
Thus $\Psi^{\times}(a)=\prod b_j$, with $b_j=\Nm_{L_j/K}(\sigma_j^{
-1}(a))$.  Let
\[F_j(\tau )=\prod_{\phi\in\Phi_j}w_{\tau\phi}^{-1}\tau w_{\phi}\quad
(\text {mod}\Aut(\mathbb{C}/K^{\text{\rm ab}})).\]
We begin by showing that $F_j(\tau )=r_K(b_j)$.
The basic properties of Artin's reciprocity law show that
\[\begin{diagram}
{\mathbb{A}}_{f,E}^{\times}&\rInto&{\mathbb{A}}_{f,\sigma L_j}^{
\times}&\rTo{\sigma_j^{-1}}&{\mathbb{A}}_{f,L_j}^{\times}&\rTo{\Nm_{L_j/K}}&
{\mathbb{A}}_{f,K}^{\times}\\
\dTo{}{r_E}&&\dTo{}{r_{\sigma L_j}}&&\dTo{}{r_{L_j}}&&\dTo{}{r_K}\\
\Gal(E^{\text{ab}}/E)&\rTo{V_{\sigma_jL_j/E}}&\sigma_
j\Gal(L_j^{\text{ab}}/L_j)\sigma_j^{-1}&\rTo{\ad\sigma_
j^{-1}}&\Gal(L_j^{\text{ab}}/L_j)&\rTo{\text{restriction}}&\Gal(K^{\text{ab}}/K)
\end{diagram}
\]
commutes.  Therefore $r_K(b_j)$ is the image of $r_E(a)$ by the three maps in the
bottom row of the diagram.  Consider $\{t_{\phi}\mid t_{\phi}=w_{
\phi}\sigma_j^{-1},\quad\phi\in\Phi_j\}$; this is a set
of coset representatives for $\sigma_j\Aut(\mathbb{C}/L_j)\sigma_j^{
-1}$ in $\Aut(\mathbb{C}/E)$, and so
$F_j(\tau )=\prod_{\phi\in\Phi_j}\sigma_j^{-1}t_{\tau\phi}^{-1}\tau
t_{\phi}\sigma_j=\sigma_j^{-1}V(\tau )\sigma_j\mod\Aut(\mathbb{C}/K^{\text{\rm ab}}
)$.
Thus $r_K(\Psi^{\times}(a))=\prod r_K(b_j)=\prod F_j(\tau )=F_{\Phi}
(\tau )$.  Clearly, $\Psi^{\times}(a)\cdot\iota\Psi^{\times}(a)\in
\chi (\tau )K^{\times}$,
and so this shows that $\Psi^{\times}(a)\in f_{\Phi}(\tau )$.
\end{proof}

\section{Start of the Proof; Tate's Result}

We shall work with the statement (1.3) rather than (1.1). The variety $ \tau
A$ has type $(K,\tau\Phi)$ because $\tau\Phi$ describes the action of $ K$ on
the tangent space to $\tau A$ at zero. Choose any $K$-linear isomorphism
$\alpha: H_{1}(A,\mathbb{Q})\to H_{1}(\tau A,\mathbb{Q})$. Then
\[
V_{f}(A)\overset{\tau}{\to}V_{f}(\tau A)\overset{(\alpha\otimes1)^{-1}}{\to
}V_{f}(A)
\]
is an $\mathbb{A}_{f,K}$-linear isomorphism, and hence is multiplication by
some $ g\in\mathbb{A}_{f,K}^{\times}$; thus $(\alpha\otimes1)\circ g=\tau$.

\begin{lemma}
For this $g$, we have
\[(\alpha\psi )(\frac {\chi (\tau )}{g\cdot\iota g}x,y)=(\tau\psi
)(x,y),\quad\text {\rm all }x,y\in V_f(\tau A).\]
\end{lemma}
\begin{proof}
By definition,
\begin{eqnarray*}
(\tau\psi )(\tau x,\tau y)&=&\tau (\psi (x,y))\qquad x,y\in V_f(A
)\\
(\alpha\psi )(\alpha x,\alpha y)&=&\psi (x,y)\qquad x,y\in V_f(A).\end{eqnarray*}
On replacing $x$ and $y$ by $gx$ and $gy$ in the second inequality, we find that
\[(\alpha\psi )(\tau x,\tau y)=\psi (gx,gy)=\psi ((g\cdot\iota g)
x,y).\]
As $\tau (\psi (x,y))=\chi (\tau )\psi (x,y)=\psi (\chi (\tau )x,
y)$, the lemma is now obvious.
\end{proof}

\begin{remark}
\begin{enumerate}
\item\ On replacing $x$ and $y$ with $\alpha x$ and $\alpha y$ in (3.1), we obtain the formula
\[\psi (\frac {\chi (\tau )}{g\cdot\iota g}x,y)=(\tau\psi )(\alpha
x,\alpha y).\]
\item On taking $x,y\in H_1(A,\mathbb{Q})$ in (3.1), we can deduce that $
\chi (\tau )/g\cdot\iota g\in K^{\times}$;
therefore $g\cdot\iota g\equiv\chi (\tau )\mod K^{\times}$.
\end{enumerate}
\end{remark}

The only choice involved in the definition of $g$ is that of $\alpha$, and $
\alpha$ is determined up to multiplication by an element of $K^{\times}$. Thus
the class of $ g$ in $\mathbb{A}_{f,K}^{\times}/K^{\times}$ depends only on
$A$ and $ \tau$. In fact, it depends only on $(K,\Phi)$ and $\tau$, because
any other abelian variety of type $(K,\Phi)$ is isogenous to $ A$ and leads to
the same class $gK^{\times}$. We define $g_{\Phi}(\tau)=gK^{\times}%
\in\mathbb{A}_{f,K}^{\times}/K^{\times}$.

\begin{proposition}
The maps $g_{\Phi}:\Aut(\mathbb{C})\to\mathbb{A}_{f,K}^{\times}/K^{\times}$ have the following properties:
\begin{enumerate}
\item\ $g_{\Phi}(\sigma\tau )=g_{\tau\Phi}(\sigma )\cdot g_{\Phi}
(\tau )$;
\item\ $g_{\Phi (\tau^{-1}|K)}(\sigma )=\tau g_{\Phi}(\sigma )$ if $
\tau K=K$;
\item\ $g_{\Phi}(\iota )=1$;
\item\ $g_{\Phi}(\tau )\cdot\iota g_{\Phi}(\tau )=\chi (\tau )K^{
\times}$.
\end{enumerate}
\end{proposition}
\begin{proof}
(a)  Choose $K$-linear isomorphisms $\alpha :H_1(A,\mathbb{Q})\to H_
1(\tau A,\mathbb{Q})$ and
$\beta :H_1(\tau A,\mathbb{Q})\to H_1(\sigma\tau A,\mathbb{Q})$, and let $
g=(\alpha\otimes 1)^{-1}\circ\tau$ and $g_{\tau}=(\beta\otimes 1)^{
-1}\circ\sigma$ so that $g$
and $g_{\tau}$ represent $g_{\Phi}(\tau )$ and $g_{\tau\phi}(\sigma
)$ respectively.  Then
\[(\beta\alpha )\otimes 1\circ (g_{\tau}g)=(\beta\otimes 1)\circ
g_{\tau}\circ (\alpha\otimes 1)\circ g=\sigma\tau ,\]
which shows that $g_{\tau}g$ represents $g_{\Phi}(\sigma\tau )$.
(b) If $(A,K\hookrightarrow\End(A)\otimes \mathbb{Q})$ has type $(K
,\Phi )$, then $(A,K\stackrel {\tau^{-1}}{\to}K\to\End(A)\otimes
\mathbb{Q})$ has type
$(K,\Phi\tau^{-1})$.  The formula in (b) can be proved by transport of structure.
(c) Complex conjugation $\iota :A\to\iota A$ is a homeomorphism (relative to the complex
topology) and so induces a $K$-linear isomorphism $\iota_1:H_1(A,
\mathbb{Q})\to H_1(A,\mathbb{Q})$.  The
map $\iota_1\otimes 1:V_f(A)\to V_f(\iota A)$ is $\iota$ again, and so on taking $
\alpha =\iota_1$, we find that $g=1$.
(d) This is proved in (3.2).
\end{proof}

Theorem 1.3 (hence also 1.1) becomes true if $f_{\Phi}$ is replaced by $
g_{\Phi}$. Our task is to show that $f_{\Phi}=g_{\Phi}$. To this end we set
\[
e_{\Phi}(\tau)=g_{\Phi}(\tau)/f_{\Phi}(\tau)\in\mathbb{A}_{f,K}^{ \times
}/K^{\times}.
\]

\begin{proposition}
The maps $e_{\Phi}:\Aut(\mathbb{C})\to\mathbb{A}_{f,K}^{\times}/K^{\times}$ have the following properties:
\begin{enumerate}
\item\ ${\rm e}_{\Phi}(\sigma\tau )={\rm e}_{\tau\Phi}(\sigma )\cdot
{\rm e}_{\Phi}(\tau )$;
\item\ ${\rm e}_{\Phi (\tau^{-1}|K)}(\sigma )=\tau {\rm e}_{\Phi}
(\sigma )$ if $\tau K=K$;
\item\ ${\rm e}_{\Phi}(\iota )=1$;
\item\ ${\rm e}_{\Phi}(\tau )\cdot\iota {\rm e}_{\Phi}(\tau )=1$;
\item\ $e_{\Phi}(\tau )=1$ if $\tau\Phi =\Phi$.
\end{enumerate}
\end{proposition}
\begin{proof}
Statements (a), (b), and (c) follow from (a), (b), and (c) of (2.6) and (3.3), and
(d) follows from (3.3d) and (2.3b).  The condition $\tau\Phi =\Phi$ in (e) means that $
\tau$
fixes the reflex field of $(K,\Phi )$ and, as we observed in \S 2, the theorems are
known to hold in that case, which means that $f_{\Phi}(\tau )=g_{
\Phi}(\tau )$.
\end{proof}

\begin{proposition}
Let $K_0$ be the maximal real subfield of $K$; then $e_{\Phi}(\tau
)\in\mathbb{A}_{f,K_0}^{\times}/K_0^{\times}$ and
$e_{\Phi}(\tau )^2=1$; moreover, $e_{\Phi}(\tau )$ depends only on the effect of $
\tau$ on $K_0$, and is $1$
if $\tau |K_0=\id$.
\end{proposition}
\begin{proof}
Replacing $\tau$ by $\sigma^{-1}\tau$ in (a), we find using (e) that $
e_{\Phi}(\tau )=e_{\Phi}(\sigma )$ if $\tau\Phi =\sigma\Phi$,
i.e., $e_{\Phi}(\tau )$ depends only on the restriction of $\tau$ to the reflex field of $
(K,\Phi )$.
>From (b) with $\tau =\iota$, we find using $\iota\Phi =\Phi\iota$ that $
e_{\iota\Phi}(\sigma )=\iota e_{\Phi}(\sigma )$.  Putting
$\tau =\iota$, then $\sigma =\iota$, in (a), we find that $e_{\Phi}
(\sigma\iota )=\iota e_{\Phi}(\sigma )$ and $e_{\Phi}(\iota\tau )
=e_{\Phi}(\tau )$.  Since
$\iota\tau$ and $\tau\iota$ have the same effect on $E$, we conclude $
e_{\Phi}(\tau )=\iota e_{\Phi}(\tau )$.  Thus
$e_{\Phi}(\tau )\in (\mathbb{A}_{f,K}^{\times}/K^{\times})^{\langle
\iota\rangle}=\mathbb{A}_{f,K_0}^{\times}/K_0^{\times}$, where $\langle
\iota\rangle =\Gal(K/K_0)$, and (d) shows that
$e_{\Phi}(\tau )^2=1$.
\end{proof}

\begin{corollary}
Part (a) of (1.1) is true; part (b) of (1.1) becomes true when $f$ is replaced by
$ef$ with $e\in\mathbb{A}_{f,K_0}^{\times}$, $e^2=1$.
\end{corollary}
\begin{proof}
Let $e\in e_{\Phi}(\tau )$.  Then $e^2\in K_0^{\times}$ and, since an element of $
K_0^{\times}$ that is a square
locally at all finite primes is a square, we can correct $e$ to achieve $
e^2=1$.
Now (1.1) is true with $f$ replaced by $ef$, but $e$ (being a unit) does not affect
part (a) of (1.1).
\end{proof}

We can now sketch the proof of the Theorems 1.1 and 1.3 --- for this, we must
prove $e_{\Phi}(\tau)=1$ for all $\tau$. It seems to be essential to prove
this simultaneously for all abelian varieties. To do this, one needs to define
a universal $e$, giving rise to all the $e_{\Phi}$. The universal $e$ is a map
into the Serre group. In $\S 4$ we review some of the theory concerning the
Serre group, and in (5.1) we state the existence of $e$. The proof of (5.1),
which requires Deligne's result (Deligne 1982a) on Hodge cycles on abelian
varieties, is carried out in \S 7 and \S 8 . The remaining step, proving that
$ e=1$, is less difficult, and is carried out in \S 6.

\section{The Serre Group}

Let $E$ be a CM-field. The \emph{Serre group\/} corresponding to $ E$ is a
pair $(S^{E},\mu^{E})$ comprising a $\mathbb{Q}$-rational torus $S^{E}$ and a
cocharacter $ \mu^{E}\in X_{*}(S^{E})$ defined over $E$ whose weight
$w^{E}\overset{\text{\textrm{df}}}{=}-(\iota+1)\mu^{E}$ is defined over
$\mathbb{Q}$. It is characterized by having the following universal property:
for any $ \mathbb{Q}$-rational torus $T$ and $\mu\in X_{*}(T)$ defined over
$E$ whose weight is $ \mathbb{Q}$-rational, there is a unique $\mathbb{Q}%
$-rational homomorphism $\rho_{\mu}:S^{E}\to T$ such that $ \rho_{\mu}\circ
\mu^{E}=\mu$.

For $\rho\in\Hom(E,\mathbb{C})$, let $[\rho]$ be the character of the torus $
E^{\times}$ defined by $\rho$. Then $\{[\rho]\mid\rho\in\Hom(E,\mathbb{C})\}$
is a basis for $X^{ *}(E^{\times})$, and $S^{E}$ is the quotient of the torus
$E^{\times}$ with
\begin{align*}
X^{*}(S^{E}) & =\{\chi\in X^{*}(E^{\times})\mid(\tau-1)(\iota+1)\chi
=0\text{\textrm{, all }}\tau\in\Aut(\mathbb{C})\}\\
X^{*}(\mu^{E}) & =\sum n_{\rho}[\rho]\mapsto n_{1}:X^{*}(S^{E})\to\mathbb{Z}%
\end{align*}
because this pair has the universal property dual to that of $(S^{E},\mu^{E}%
)$. In particular, there is a canonical homomorphism $E^{\times}\to S^{E}$,
and it is known (cf. Serre 1968, II) that the kernel of the map is the Zariski
closure of any sufficiently small subgroup $U$ of finite index in
$\mathcal{O}_{E}^{\times}$.

When $E$ is Galois over $\mathbb{Q}$, the action of $\sigma\in\Gal(
E/\mathbb{Q})$ on $E$ defines an automorphism $\tilde{\sigma}$ of the torus
$S^{E}$, whose action on characters is
\[
\sum n_{\rho}[\rho]\mapsto\sum n_{\rho}[\rho\sigma]=\sum n_{\rho\sigma^{-1}%
}[\rho].
\]

\begin{lemma}
Let $E_0$ be the maximal real subfield of $E$; there is an exact sequence of
algebraic tori
\[\begin{CD}
1@>>>E_0^{\times}@>{\left(\begin{smallmatrix}\text{incl.}\\
\Nm_{E_0/{\mathbb{Q}}}\end{smallmatrix}\right)}>>E^{\times}\times{\mathbb{Q}}^{\times}
@>{\left(\text{can.},w^E\right)}>>S^E\to1.
\end{CD}\]
\end{lemma}
\begin{proof}
It suffices to show that the sequence becomes exact after the functor $
X^{*}$
has been applied.  As
\begin{eqnarray*}
X^{*}(E_0)&=&\{\sum n_{\rho'}[\rho']\mid\rho'\in\Hom(E_0,\mathbb{C}
)\}\\
X^{*}(E^{\times}\times \mathbb{Q}^{\times})&=&\{\sum n_{\rho}[\rho
]+n\mid\rho\in\Hom(E,\mathbb{C})\}\\
X^{*}(S^E)&=&\{\sum n_{\rho}[\rho ]\mid n_{\rho}+n_{\iota\rho}=\text{\rm constant}
\}\\
X^{*}(\left(\begin{array}{c}
\text{incl.}\\
\Nm_{E_0/\mathbb{Q}}\end{array}
\right))&=&\sum n_{\rho}[\rho ]+n\mapsto\sum n_{\rho}[\rho |E_0]+n
\sum_{\rho'}[\rho']\\
X^{*}((\text{\rm can}.,w^E))&=&\sum n_{\rho}[\rho ]\mapsto\sum n_{
\rho}[\rho ]-(n_1+n_{\iota})\end{eqnarray*}
this is trivial.
\end{proof}

\begin{lemma}
The map $\Nm_{E/\mathbb{Q}}:E^{\times}\to \mathbb{Q}^{\times}$ factors through $
S^E$, and gives rise to a
commutative diagram
\[\begin{diagram}[heads=LaTeX]
S^E&&\rTo{1+\iota}&&S^E\\
&\rdTo{}{\Nm_{E/{\mathbb{Q}}}}&&\ruTo{}{-w^E}\\
&&{\mathbb{Q}}^{\times}.&\end{diagram}\]
\end{lemma}
\begin{proof}
The map $X^{*}(\Nm_{E/\mathbb{Q}})$ is $n\mapsto n\sum [\rho ]$, which clearly factors through
$X^{*}(S^E)\subset X^{*}(E^{\times})$.  Moreover, the endomorphisms
\[X^{*}(-w^E\circ\Nm_{E/\mathbb{Q}})=(\sum n_{\rho}[\rho ]\mapsto n_
1+n_{\iota}\longmapsto (n_1+n_{\iota})\sum n_{\rho}[\rho ])\]
\[X^{*}(1+\iota )=(\sum n_{\rho}[\rho ]\mapsto\sum n_{\rho})([\rho
]+[\iota\rho ])\mapsto\sum (n_{\rho}+n_{\iota\rho})[\rho ]=(n_1+n_{
\iota})\sum [\rho ]\]
are equal.
\end{proof}

Let $E_{1}\supset E_{2}$ be CM-fields. The norm map $E_{1}^{\times}\to
E_{2}^{\times}$ induces a norm map $\Nm_{E_{1}/E_{2}}:S^{E_{1}}\to S^{E_{2}}$
which is the unique $\mathbb{Q}$-rational homomorphism such that
$\Nm_{E_{1}/E_{2}}\circ\mu^{E_{1}}=\mu^{E_{2}}$. The following diagram commutes:%

\[
\begin{CD}
1@>>>(E_1)_0^{\times}@>>>E_1^{\times}\times{\mathbb{Q}}^{\times}@>>>S^{E_1}@>>>1\\
@.@VV{\Nm}V@VV{\Nm\times\id}V@VV{\Nm}V\\
1@>>>(E_2)_0^{\times}@>>>E_2^{\times}\times{\mathbb{Q}}^{\times}@>>>S^{E_2}@>>>1.
\end{CD}
\]

\begin{remark}
The Serre group can be defined for all fields of finite degree over $
\mathbb{Q}$.  If $L$
contains a CM-field and $E$ is the maximal such subfield, then $\Nm_{
L/K}:S^L\stackrel {\approx}{\to}S^E$;
if not, then $\Nm_{L/\mathbb{Q}}:S^L\stackrel {\approx}{\to}S^{\mathbb{Q}}
=\mathbb{Q}^{\times}$.
\end{remark}

Let $(K,\Phi)$ be a CM-type with $K\subset\mathbb{C}$. Write $T=\Res_{
K/\mathbb{Q}}(\mathbb{G}_{m})$, and define $\mu_{\Phi}\in X_{*}(T)$ by the
condition
\[
[\rho]\circ\mu_{\Phi}=\left\{
\begin{array}
[c]{rll}%
\id & ,\quad\rho\in\Phi & \\
1, & \quad\rho\notin\Phi. &
\end{array}
\right.
\]
Thus, $\mu_{\Phi}$ is the map
\[
\begin{diagram}
{\mathbb{C}}^{\times}&\to&T({\mathbb{C}})=(K\otimes_{{\mathbb{Q}}}{\mathbb{C}})^{
\times}&=&\prod_{\phi\in\Phi}{\mathbb{C}}^{\times}\times\prod_{\phi\notin
\Phi}{\mathbb{C}}^{\times}\\
z&&\mapsto&&(z,\ldots ,z,1,\ldots ,1).\end{diagram}
\]
The weight of $\mu_{\Phi}$ is the map induced by $x\mapsto x^{-1}
:\mathbb{Q}^{\times}\hookrightarrow K^{\times}$, which is defined over
$\mathbb{Q}$, and $\mu_{\Phi}$ itself is defined over the reflex field of $
(K,\Phi)$. There is therefore, for any CM-field $E$ containing the reflex
field of $( K,\Phi)$, a unique $\mathbb{Q}$-rational homomorphism $\rho_{\Phi
}:S^{E}\to T$ such that $ \mu_{\Phi}=\rho_{\Phi}\circ\mu^{E}$. From now on, we
assume $E$ to be Galois over $\mathbb{Q}$.

\begin{lemma}
\begin{enumerate}
\item\ $\tau\mu_{\Phi}=\mu_{\tau\Phi}$, $\tau\in\Aut(\mathbb{C})$.
\item\ Let $\tau\in\Aut(\mathbb{C})$ be such that $\tau K=K$, so that $
\tau$ induces an
automorphism $\tilde{\tau}$ of $T$; then $\tilde{\tau}\circ\mu_{\Phi}
=\mu_{\Phi\tau^{-1}}$.
\end{enumerate}
\end{lemma}
\begin{proof}
(a) Consider the canonical pairing
\[\langle\cdot ,\cdot\rangle :X^{*}(T)\times X_{*}(T)\to \mathbb{Z}
.\]
By definition, for $\rho\in\Hom(K,\mathbb{C})$,
\[\langle [\rho ],\mu_{\Phi}\rangle =\left\{\begin{array}{ll}
1&\text{\rm if }\rho\in\Phi\\
0&\text{\rm otherwise}.\end{array}
\right.\]
For $\tau\in\Aut(\mathbb{C})$,
\[\langle [\rho ],\tau\mu_{\Phi}\rangle =\langle\tau^{-1}[\rho ],
\mu_{\Phi}\rangle =\langle [\tau^{-1}\rho ],\mu_{\Phi}\rangle ,\]
which equals $\langle [\rho ],\mu_{\tau\Phi}\rangle .$
(b)
\[[\rho ]\circ\tilde{\tau}\circ\mu_{\Phi}=[\rho\tau ]\circ\mu_{\Phi}
=\left\{\begin{array}{rcl}
\id&\text{\rm \ if }&\rho\tau\in\Phi\\
1&\text{\rm \ if }&\rho\tau\notin\Phi .\end{array}
\right.\]
Thus (b) is clear.
\end{proof}

\begin{proposition}
\begin{enumerate}
\item\ For any $\tau\in\Aut(\mathbb{C})$, $\rho_{\Phi}\circ\tilde{\tau}^{
-1}=\rho_{\tau\Phi}$.
\item\ If $\tau K=K$, then $\tilde{\tau}\circ\rho_{\Phi}=\rho_{\Phi
\tau^{-1}}$.
\end{enumerate}
\end{proposition}
\begin{proof}
(a) We shall show that $\tilde{\tau}^{-1}\circ\mu^E=\tau (\mu^E)$; from this it follows that
\[\begin{array}{rclll}
\rho_{\Phi}\circ\tilde{\tau}^{-1}\circ\mu^E&=&\rho_{\Phi}\circ (\tau
\mu^E)\\
&=&\tau (\rho_{\Phi}\circ\mu^E)&\qquad&(\rho_{\Phi}\text{\rm \ is $
\mathbb{Q}$-rational})\\
&=&\tau (\mu_{\Phi})&\qquad&\text{\rm (definition of $\rho_{\Phi}$)}\\
&=&\mu_{\tau\Phi}&&\text{\rm (4.4a),}\end{array}
\]
which implies that $\rho_{\Phi}\circ\tau^{-1}=\rho_{\tau\Phi}$.  It remains to show that
$X^{*}(\tilde{\tau}^{-1}\circ\mu^E)=X^{*}(\tau\mu^E)$, but
\[X^{*}(\tilde{\tau}^{-1}\circ\mu^E)=X^{*}(\mu^E)\circ X^{*}(\tilde{
\tau}^{-1})=(\sum n_{\rho}[\rho ]\mapsto\sum n_{\rho}[\rho\tau^{-
1}]\mapsto n_{\tau})\]
and
\[X^{*}(\tau\mu^E)=\tau (X^{*}(\mu^E))=(\sum n_{\rho}[\rho ]\mapsto
\sum n_{\rho}[\tau^{-1}\rho ]\mapsto n_{\tau}).\]
(b)
\[\begin{array}{ccccc}
\tilde{\tau}\circ\rho_{\Phi}\circ\mu^E&=&\tilde{\tau}\circ\mu_{\Phi}&
\quad&\text{\rm (definition of $\rho_{\Phi})$}\\
&=&\mu_{\Phi\tau^{-1}}&\qquad&\text{\rm (4.4b)}.\end{array}
\]
\end{proof}

\begin{proposition}
For $E_1\supset E_2$,
\[\begin{diagram}[heads=LaTeX,size=3em]
S^{E_1}\\
\dTo{\Nm_{E_1/E_2}}&\rdTo{\rho_{\Phi}}\\
S^{E_2}&\rTo{\rho_{\Phi}}&T
\end{diagram}\]
commutes.
\end{proposition}
\begin{proof}
We have
\[((\rho_{\Phi})_2\circ\Nm_{E_1/E_2})\circ\mu^{E_1}=(\rho_{\Phi})_
2\circ\mu^{E_2}=\mu_{\Phi}.\]
\end{proof}

\begin{proposition}
Let $E$ be a CM-field, Galois over $\mathbb{Q}$, and consider all maps $
\rho_{\Phi}$ for $\Phi$ running
through the CM-types on $E$; then $\cap\Ker(\rho_{\Phi})=1$.
\end{proposition}
\begin{proof}
We have to show that $\sum\text{\rm Im}(X^{*}(\rho_{\Phi}))=X^{*}
(S^E)$; but the left hand side contains
$\sum_{\Phi}[\phi ]$ for all CM-types on $E$, and these elements generate $
X^{*}(S^E)$.
\end{proof}

\begin{proposition}
Let $K_1\supset K_2$ be CM-fields,  and let $\Phi_1$ and $\Phi_2$ be CM-types for $
K_1$ and $K_2$
respectively such that $\Phi_1|K_2=\Phi_2$.  Then, for any CM-field $
E$ containing the
reflex field of $(K_1,\Phi_1)$, the composite of
\[S^E\stackrel {\rho_{\Phi_2}}{\to}K_2^{\times}\hookrightarrow K_
1^{\times}\]
is $\rho_{\Phi_1}$.
\end{proposition}
\begin{proof}
Let $i:K_2^{\times}\hookrightarrow K_1^{\times}$ be the inclusion map.  Then $
i\circ\mu_{\Phi_2}=\mu_{\Phi_1}$ and so
$i\circ\rho_{\Phi_2}\circ\mu^E=i\circ\mu_{\Phi_2}=\mu_{\Phi_1}$, which shows that $
i\circ\rho_{\Phi_2}=\rho_{\Phi_1}$.
\end{proof}

\section{Definition of $e^{E}$}

\begin{proposition}
Let $E\subset \mathbb{C}$ be a CM-field, Galois over $\mathbb{Q}$.  Then there exists a unique map
$e^E:\Aut(\mathbb{C})\to S^E(\mathbb{A}_f)/S^E(\mathbb{Q})$ such that, for all CM-types $
(K,\Phi )$ whose reflex
fields are contained in $E$, $e_{\Phi}(\tau )=\rho_{\Phi}(e(\tau
))$.
\end{proposition}
\begin{proof}
The existence of $e^E$ will be shown in \S 7  and \S 8.  The uniqueness follows
from (4.7) for this shows that there is an injection $S^E\rInto{
(\rho_{\Phi})}\prod T_{\Phi}$ where
$T_{\Phi}=\Res_{E/\mathbb{Q}}\mathbb{G}_m$ and the product is over all CM-types on $
E$.  Thus
\[S^E(\mathbb{A}_f)/S^E(\mathbb{Q})\hookrightarrow\prod T_{\Phi}(\mathbb{A}_
f)/T_{\Phi}(\mathbb{Q})=\prod \mathbb{A}_{f,E}^{\times}/E^{\times},\]
and so any element $a\in S^E(\mathbb{A}_f)/S^E(\mathbb{Q})$ is determined by the set $
(\rho_{\Phi}(a))$.
\end{proof}

\begin{proposition}
The family of maps $e^E:\Aut(\mathbb{C})\to S^E(\mathbb{A}_f)/S^E(\mathbb{Q}
)$ has the following properties:
\begin{enumerate}
\item\ $e^E(\sigma\tau )=\tilde{\tau}^{-1}e^E(\sigma )\cdot e^E(\tau
),\quad\sigma ,\tau\in\Aut(\mathbb{C})$;
\item\ if $E_1\supset E_2$, then
\[\begin{diagram}[heads=LaTeX,size=3em]
\Aut({\mathbb{C}})&\rTo{e^{E_1}}&S^{E_1}({\mathbb{A}}_f)/S^{E_1}({\mathbb{Q}})\\
&\rdTo{e^{E_2}}&\dTo>{\Nm}\\
&&S^{E_2}({\mathbb{A}}_f)/S^{E_2}({\mathbb{Q}})
\end{diagram}\]
commutes.
\item\ $e^E(\iota )=1$;
\item\ $e(\tau )\cdot\tilde{\iota }e(\tau )=1,\quad\tau\in\Aut(
\mathbb{C})$;
\item\ $e^E|\Aut(\mathbb{C}/E)=1$.
\end{enumerate}
\end{proposition}
\begin{proof}
(a) We have to check that $\rho_{\Phi}(e(\sigma\tau ))=\rho_{\Phi}
(\tilde{\tau}^{-1}e(\sigma )\cdot e(\tau ))$ for all $(K,\Phi )$.  But
$\rho_{\Phi}(e^E(\sigma\tau ))=e_{\Phi}(\sigma\tau )$ and
\[\rho_{\Phi}(\tilde{\tau}^{-1}e^E(\sigma )e^E(\tau ))=\rho_{\Phi}
(\tilde{\tau}^{-1}e^E(\sigma ))\rho_{\Phi}(e^E(\tau ))=\rho_{\tau
\Phi}(e^E(\sigma ))\cdot\rho_{\Phi}(e^E(\tau ))=e_{\tau\Phi}(\sigma
)e_{\Phi}(\tau );\]
thus the equality follows from (3.4a).
(b) This follows from (4.6) and the definition of $e^E$.
(c) $\rho_{\Phi}(e^E(\iota ))=e_{\Phi}(\iota )=1$ (by 3.4c), and so $
e^E(\iota )=1$.
(d) This follows from (3.4d).
(e) Assume $\tau$ fixes $E$; then $\tau\Phi =\Phi$ whenever $E$ contains the reflex field of
$(K,\Phi )$, and so $\rho_{\Phi}(e^E(\tau ))=e_{\Phi}(\tau )=1$  by (3.4e).
\end{proof}

\begin{remark}
\begin{enumerate}
\item\ Define $\varepsilon^E(\tau )=e(\tau^{-1})^{-1}$; then the maps $
\varepsilon^E$ satisfy the same conditions
(b), (c), (d), and (e) of (5.2) as $e^E$, but (a) becomes the condition
$\varepsilon^E(\sigma\tau )=\tilde{\sigma}\varepsilon^E(\tau )\cdot
\varepsilon^E(\sigma )$: $\varepsilon^E$ is a crossed homomorphism.
\item\ Condition (b) shows that $e^E$ determines $e^{E'}$ for all $
E'\subset E$.  We extend
the definition of $e^E$ to all CM-fields $E\subset \mathbb{C}$ by letting $
e^E=\Nm_{E_1/E}\circ e^{E_1}$
for any Galois CM-field $E_1$ containing $E$.
\item\ Part (d) of (5.2) follows from the remaining parts, as is clear from
the following diagram:
\[\begin{diagram}[heads=LaTeX,size=3em]
\Aut({\mathbb{C}})&\rTo{e^E}&S^E({\mathbb{A}}_f)/S^E({\mathbb{Q}})&\rTo{1+\iota}&S^E({\mathbb{A}}_f)/S^E({\mathbb{Q}})\\
&\rdTo{e^{{\mathbb{Q}}[i]}=1}&\dTo{}{\Nm_{E/{\mathbb{Q}}[i]}}&\rdTo{\Nm_{E/{\mathbb{Q}}}}&\uTo{}{-w^E}\\
&&S^{{\mathbb{Q}}[i]}({\mathbb{A}}_f)/S^{{\mathbb{Q}}[i]}({\mathbb{Q}})&\rTo{\Nm_{{\mathbb{Q}}
[i]/{\mathbb{Q}}}}&S^{{\mathbb{Q}}}({\mathbb{A}}_f)/S^{{\mathbb{Q}}}({\mathbb{Q}})
\end{diagram}\]
The right hand triangle is (4.2).  (We can assume $E\supset \mathbb{Q}
[i]$; $S^{\mathbb{Q}[i]}=\mathbb{Q}[i]^{\times}$,
$S^{\mathbb{Q}}=\mathbb{Q}^{\times}$).  In his (original) letter to Langlands (see Deligne 1979),  Deligne
showed that the difference between the motivic Galois group and the
Taniyama group was measured by a family of crossed homomorphisms $
(e^E)$
having properties (b), (c), and (e) of (5.2).  After seeing Tate's result he
used the above diagram to show that his maps $e^E$ had the same properties
as Tate's $e_{\Phi}(\tau )$, namely, $e^E(\tau )\cdot\iota e^E(\tau
)=1$, $e^E(\tau )^2=1$.
\end{enumerate}
\end{remark}

\section{Proof that $e=1$}

We replace $e$ with $\tau\mapsto e(\tau^{-1})^{-1}$.

\begin{proposition}
Suppose there are given crossed homomorphisms $e^E:\Aut(\mathbb{C})
\to S^E(\mathbb{A}_f)/S^E(\mathbb{Q})$,
one for each CM-field $E\subset \mathbb{C}$, such that
\begin{enumerate}
\item\ $e^E(\iota )=1$, all $E$;
\item\ $e^E|\Aut(\mathbb{C}/E)=1$;
\item\ if $E_1\supset E_2$ then
\[\begin{diagram}[heads=LaTeX,size=3em]
\Aut({\mathbb{C}}/{\mathbb{Q}})&\rTo{e^{E_1}}&S^{E_1}({\mathbb{A}}_f)/S^{E_1}({\mathbb{Q}})\\
&\rdTo{e^{E_2}}&\dTo{}{\Nm_{E_1/E_2}}\\
&&S^{E_2}({\mathbb{A}}_f)/S^{E_2}({\mathbb{Q}})&
\end{diagram}\]
commutes.
\end{enumerate}
Then $e^E=1$ --- i.e., $e^E(\tau )=1$ for all $\tau$ --- for all $
E$.
\end{proposition}
\begin{proof}
Clearly, it suffices to show that $e^E=1$ for all sufficiently large $
E$ --- in
particular, for those that are Galois over $\mathbb{Q}$.
The crossed homomorphism condition is that
\[e(\sigma\tau )=\tilde{\sigma }e(\tau )\cdot e(\sigma ).\]
Condition (b) implies that $e^E(\tau )=e^E(\tau')$ if $\tau |E=\tau'
|E$.  In particular,
$e^E(\iota\tau )=e^E(\tau\iota )$ for all $\tau\in\Aut(\mathbb{C})$.  Since
\[\left\{\begin{array}{rcl}
e^E(\tau\iota )&=\tau e^E(\iota )\cdot e^E(\tau )&=e^E(\tau )\\
e^E(\iota\tau )&=\iota e^E(\tau )\cdot e^E(\iota )&=\iota e^E(\tau
)\end{array}
\right.\]
we conclude that $\iota e^E(\tau )=e^E(\tau )$.
\begin{lemma}
Assume that $E$ is Galois over $\mathbb{Q}$, and let $\langle\iota\rangle$ be the subgroup of $\Gal
(E/\mathbb{Q})$
generated by $\iota |E$.
\begin{enumerate}
\item\ There is an exact commutative diagram
\[\begin{CD}
1@>>>{\mathbb{Q}}^{\times}@>>>S^E({\mathbb{Q}})^{\langle\iota\rangle}@>>>
\mu_2(E_0)@>{\Nm_{E_0/{\mathbb{Q}}}}>>\mu_2({\mathbb{Q}})\\
@.@VVV@VVV@VVV@VVV\\
1@>>>{\mathbb{A}}_f^{\times}@>>>S^E({\mathbb{A}}_f)^{\langle\iota\rangle}@>>>
\mu_2({\mathbb{A}}_{f,E_0})@>{\Nm}>>\mu_2({\mathbb{A}}_f)
\end{CD}
\]
where $\mu_2(R)$ denotes the set of square roots of $1$ in a ring $
R$.
\item\ The canonical map
\[H^1(\langle\iota\rangle ,S^E(\mathbb{Q}))\to H^1(\langle\iota\rangle
,S^E(\mathbb{A}_f))\]
is injective.
\end{enumerate}
\end{lemma}
\begin{proof}
>From (4.1) we obtain cohomology sequences
\[\begin{diagram}[width=1em]
1&\to&E_0^{\times}&\to&E_0^{\times}\times{\mathbb{Q}}^{\times}&\to&S^
E({\mathbb{Q}})^{\langle\iota\rangle}&\to&\mu_2(E_0)&\rTo{\Nm}&
\mu_2({\mathbb{Q}})&\to&H^1(\langle\iota\rangle,S^E({\mathbb{Q}}))&\to&E_
0^{\times}/E_0^{\times2}\\
&&\dTo&&\dTo&&\dTo&&\dTo&&\dTo&&\dTo&&\dTo\\
1&\to&{\mathbb{A}}_{f,E_0}^{\times}&\to&{\mathbb{A}}_{f,E_0}^{\times}\times{\mathbb{A}}_
f^{\times}&\to&S^E({\mathbb{A}}_f)^{\langle\iota\rangle}&\to&\mu_2({\mathbb{A}}_{
f,E_0})&\rTo{\Nm}&\mu_2({\mathbb{A}}_f)&\to&H^1(\langle\iota
\rangle,S^E({\mathbb{A}}_f))&\to&{\mathbb{A}}_{f,E_0}^{\times}/{\mathbb{A}}_{
f,E_0}^{\times2}\end{diagram}\]
It is easy to extract from this the diagram in (a).  For (b), let
$\gamma\in H^1(\langle\iota\rangle ,S^E(\mathbb{Q}))$ map to zero in $
H^1(\langle\iota\rangle ,S^E(\mathbb{A}_f))$.  As $E_0^{\times}/E_0^{
\times 2}\to \mathbb{A}_{f,E_0}^{\times}/\mathbb{A}_{f,E_0}^{\times 2}$ is
injective (an element of $E_0$ that is a square in $E_{0,v}$ for all finite primes is
a square in $E_0$), we see that $\gamma$ is the image of $\pm 1\in
\mu_2(\mathbb{Q})$.  The map
$\Nm_{E_0/\mathbb{Q}}:\mu_2(E_0)\to\mu_2(\mathbb{Q})$ sends $-1$ to $
(-1)^{[E_0:\mathbb{Q}]}$.  If $[E_0:\mathbb{Q}]$ is odd, it is
surjective, and therefore $\gamma =0$.  Suppose therefore that $[
E_0:\mathbb{Q}]$ is even, and
that $\gamma$ is the image of $-1$.  The assumption that $\gamma$ maps to zero in
$H^1(\langle\iota\rangle ,S^E(\mathbb{A}_f))$ then implies that $-1
\in \mathbb{Q}_{\ell}$ is in the image of $\Nm:E_0\otimes \mathbb{Q}_{
\ell}\to \mathbb{Q}_{\ell}$
for all $\ell$; but this is impossible, since for some $\ell$, $[
E_{0v}:\mathbb{Q}_{\ell}]$ will be even for
one (hence all) $v$ dividing $\ell$.
\end{proof}
Part (b) of the lemma shows that
\[S^E(\mathbb{A}_f)^{\langle\iota\rangle}/S^E(\mathbb{Q})^{\langle\iota
\rangle}=(S^E(\mathbb{A}_f)/S^E(\mathbb{Q}))^{\langle\iota\rangle}.\]
The condition $\iota e^E(\tau )=e^E(\tau )$ shows that $e^E$ maps into the right hand group,
and we shall henceforth regard it as mapping into the left hand group.
>From part (a) we can extract an exact sequence
\[1\to \mathbb{A}_f^{\times}/\mathbb{Q}^{\times}\stackrel w{\to}S^E(\mathbb{A}_
f)^{\langle\iota\rangle}/S^E(\mathbb{Q})^{\langle\iota\rangle}\to\mu_
2(\mathbb{A}_{f,E_0})/\mu_2(E_0).\]
Now assume that $E\supset \mathbb{Q}[i]$, so that $E=E_0[i]$.  We show first that the image
of $e^E(\tau )$ in $\mu_2(\mathbb{A}_{f,E_0})/\mu_2(E_0)$ is $1$.  Let $
\varepsilon$ represent the image; then $\varepsilon =(\varepsilon_
v)$,
$\varepsilon_v=\pm 1$, and $\varepsilon$ itself is defined up to sign.  We shall show that, for any two
primes $v_1$ and $v_2$, $\varepsilon_{v_1}=\varepsilon_{v_2}$.  Choose a totally real quadratic extension $
E_0'$ of
$E_0$ in which $v_1$ and $v_2$ remain prime, and let $E'=E_0'[i]$.  Let $
\varepsilon'$ represent
the image of $e^{E'}(\tau )$ in $\mu_2(\mathbb{A}_{f,E_0'})/\mu_2(E_
0')$.  Then condition (c) shows that
$\Nm_{E_0'/E_0}\varepsilon'$ represents the image of $e^E(\tau )$, and so $\Nm_{
E_0'/E_0}\varepsilon'=\pm\varepsilon$.  But if
$v_i'|v_i$, then $\Nm_{E_{0,v_i'}'/E_{0,v_i}}=1$ for $i=1,2$.
It follows that $e^E$ factors through $w(\mathbb{A}_f^{\times}/\mathbb{Q}^{
\times})$.  Consider,
\[\begin{CD}
E:\qquad{1}@>>>{\mathbb{A}}_f^{\times}/{\mathbb{Q}}^{\times}@>w>>S^E({\mathbb{A}}_f)/S^E({\mathbb{Q}})\\
@.@VV{\id}V@VVV\\
{\mathbb{Q}}[i]:\qquad{1}@>>>{\mathbb{A}}_f^{\times}/{\mathbb{Q}}^{\times}@>w>>S^{{\mathbb{Q}}[i]}({\mathbb{A}}_f)/S^{{\mathbb{Q}}[i]}({\mathbb{Q}}).
\end{CD}\]
According to (c), $e^E(\tau )$ maps to $e^{\mathbb{Q}[i]}(\tau )$ under the right hand arrow, which
according to (a) and (b), is $1$.  As $e^E(\tau )$ lies in $\mathbb{A}_
f^{\times}/\mathbb{Q}^{\times}$, and the map from
there into $S^{\mathbb{Q}[i]}(\mathbb{A}_f)/S^{\mathbb{Q}[i]}(\mathbb{Q})$ is injective, this shows that
$e^E(\tau )=1$.
\end{proof}

\begin{remark}
The argument used in the penultimate paragraph of the above proof is that
used by Shih 1976, p101, to complete his proof of his special case of (1.1).
For the argument in the final paragraph, cf. 5.3c.  These two arguments
were all that was lacking in the original version Deligne 1979b of Deligne
1982.
\end{remark}

\section{Definition of $f^{E}$}

We begin the proof of (5.1) by showing that there is a universal $ f$, giving
rise to the $f_{\Phi}$.

Let $E\subset\mathbb{C}$. The Weil group $W_{E/\mathbb{Q}}$ of $E/\mathbb{Q}$
fits into an exact commutative diagram:
\[
\begin{diagram}
1&\rTo&{\mathbb{A}}_{E}^{\times}/E^{\times}&\rTo&W_{E/{\mathbb{Q}}}&\rTo&\Hom
(E,{\mathbb{C}})&\rTo&1\\
&&\dOnto{}{\text{rec}_E}&&\dOnto&&\|\\
1&\rTo&\Gal(E^{\text{ab}}/E)&\rTo&\Hom(E^{\text{ab}},{\mathbb{C}}
)&\rTo&\Hom(E,{\mathbb{C}})&\rTo&1\end{diagram}
\]
(see Tate 1979). Assume that $E$ is totally imaginary. Then $E_{\infty
}^{\times}E^{\times}\subset\Ker(\rec_{E})$, and so we can divide out by this
group and its image in $W_{E/\mathbb{Q}}$ to obtain the exact commutative
diagram:
\[
\begin{diagram}
1&\rTo&{\mathbb{A}}_{f,E}^{\times}/E^{\times}&\rTo&W_{f,E/{\mathbb{Q}}}&\rTo&\Hom
(E,{\mathbb{C}})&\rTo&1\\
&&\dOnto{}{r_E}&&\dOnto&&\|\\
1&\rTo&\Gal(E^{\text{ab}}/E)&\rTo&\Hom(E^{\text{ab}},{\mathbb{C}}
)&\rTo&\Hom(E,{\mathbb{C}})&\rTo&1\end{diagram}
\]
Assume now that $E$ is a CM-field Galois over $\mathbb{Q}$. The cocharacter $
\mu^{E}$ is defined over $E$, and gives rise to a map $\mu^{E}(R):R^{\times
}\to S^{E}(R)$ for any $E$-algebra $R$. Choose elements $w_{\sigma}\in
W_{E/\mathbb{Q}}^{f}$, one for each $ \sigma\in\Hom(E,\mathbb{C})$, such that
\[
w_{\sigma}|E=\sigma,\quad w_{\iota\sigma}=\tilde{\iota}w_{\sigma
}\text{\textrm{\ all }} \sigma,
\]
where $\tilde{\iota}$ maps to $\iota\in\Hom(E^{\text{\textrm{ab}}},\mathbb{C}
)$ (cf. \S 2). Let $\tau\in\Aut(\mathbb{C})$ and let $\tilde{\tau}\in
W_{f,E/\mathbb{Q}}$ map to $\tau|E^{ab}$. Then $ w_{\tau\sigma}^{-1}%
\circ\tilde{\tau}\circ w_{\sigma}\in\mathbb{A}_{ f,E}$, and we define
\[
f(\tau)=\prod_{\sigma\in\Hom(E,\mathbb{C})}(\sigma^{-1}\mu^{E})(w_{ \tau
\sigma}^{-1}\tilde{\tau}w_{\sigma})\mod S^{E}(E).
\]
Thus $f$ is a map $\Aut(\mathbb{C})\to S^{E}(\mathbb{A}_{f,E})/S^{E}(E)$.

\begin{proposition}
Let $(K,\Phi )$ be a CM-type whose reflex field is contained in $
E$, and let
$T=\Res_{K/\mathbb{Q}}\mathbb{G}_m$.  Identify $T(\mathbb{A}_f)/T(\mathbb{Q}
)$ with a subgroup of $T(\mathbb{A}_{f,E})/T(E)$.
Then
\[\rho_{\Phi}(f(\tau ))=f_{\Phi}(\tau ).\]
\end{proposition}
\begin{proof}
Because of (4.8), it suffices to show this with $K=E$.
\begin{lemma}
With the above notations,
\[f_{\Phi}(\tau )=\prod_{\phi\in\Phi}w_{\tau\phi}^{-1}\circ\tilde{
\tau}\circ w_{\phi}\mod E^{\times}.\]
\end{lemma}
\begin{proof}
Let $f'$ denote the right hand side.  Then $r_E(f')=F_{\Phi}(\tau
)$ (obviously), and the
same argument as in the proof of (2.3) shows that $f'\cdot\iota f'
=\chi (\tau )$.\end{proof}
We now assume that $E=K$, $E/\mathbb{Q}$ Galois.  Write $i$ for the map $
T(\mathbb{Q})\to T(E)$
induced by $\mathbb{Q}\hookrightarrow E$; then, for any $\rho\in\Hom
(E,\mathbb{C})$ and $a\in T(\mathbb{Q})=E^{\times}$, $[\rho ](i(a))=\rho
a$.
Thus $[\rho ](i(f_{\Phi}(\tau )))=\rho f_{\Phi}(\tau )=f_{\Phi\rho^{
-1}}(\tau )$ by (2.6b).  On the other hand,
\begin{eqnarray*}
[\rho ](\rho_{\Phi}(f(\tau )))&=&[\rho ]\prod_{\sigma}\rho_{\Phi}
\circ (\sigma^{-1}\mu^E)(w_{\tau\sigma}^{-1}\tilde{\tau }w_{\sigma}
)\\
&=&[\rho ]\prod_{\sigma}\sigma^{-1}(\rho_{\Phi}\circ\mu^E)(w_{\tau
\sigma}^{-1}\tilde{\tau }w_{\sigma})\\
&=&[\rho ](\prod_{\sigma}\sigma^{-1}\mu_{\Phi}(w_{\tau\sigma}^{-1}
\tilde{\tau }w_{\sigma}))\\
&=&\prod_{\sigma}([\rho ]\circ\mu_{\sigma^{-1}\Phi})(w_{\tau\sigma}^{
-1}\tilde{\tau }w_{\sigma})\qquad \text{( by 4.4a)}\\
&=&\prod_{\sigma\text{ such that }\rho\in\sigma^{-1}\Phi}w_{
\tau\sigma}^{-1}\tilde{\tau }w_{\sigma}\\
&=&\prod_{\sigma\in\Phi\rho^{-1}}w_{\tau\sigma}^{-1}\tilde{\tau }
w_{\sigma}\\
&=&f_{\Phi\rho^{-1}}(\tau ).\end{eqnarray*}
\end{proof}

\begin{corollary}
\begin{enumerate}
\item\ $f(\tau )$ depends only on $E$ and $\tau$; we have therefore defined maps
$f^E:\Aut(\mathbb{C})\to S^E(\mathbb{A}_{f,E})/S^E(E)$, one for each CM-field (Galois over $
\mathbb{Q}$);
\item\ $f^E(\sigma\tau )=\tilde{\tau}^{-1}f^E(\sigma )\cdot f^E(\tau
),\quad\sigma ,\tau\in\Aut(\mathbb{C})$;
\item\ if $E_1\supset E_2$, then
\[\begin{diagram}[heads=LaTeX,size=3em]
\Aut({\mathbb{C}})&\rTo{e^{E_1}}&S^{E_1}({\mathbb{A}}_f)/S^{E_
1}({\mathbb{Q}})\\
&\rdTo{}{e^{E_2}}&\dTo\\
&&S^{E_2}({\mathbb{A}}_f)/S^{E_2}({\mathbb{Q}})
\end{diagram}
\]
commutes;
\item\ $f^E(\iota )=1$;
\item\ $f^E(\tau )\cdot\tilde{\imath }f^E(\tau )=w^E(\tau )^{-1}$;
\item\ $\sigma f^E(\tau )=f^E(\tau )$  for all $\sigma\in\Gal(E/\mathbb{Q}
)$.
\end{enumerate}
\end{corollary}
\begin{proof}
(a) $f(\tau )$ is the unique element of $S^E(\mathbb{A}_{f,E})/S^E(
E)$ such that $\rho_{\Phi}(f(\tau ))=f_{\Phi}(\tau )$
for all $(K,\Phi )$.  (Cf. the proof of the uniqueness of $e^E$ in (5.1).)
(b), (c), (d), (e).  These are proved as (a), (b), (c), (d) of (5.2).
(f).$\rho_{\Phi}(\sigma f^E(\tau ))=\sigma (\rho_{\Phi}(f^E(\tau )))=\sigma f_{\Phi}^E(\tau )=f_{\Phi}^E(\tau )$.
\end{proof}

\begin{remark}
Let $\bar {w}_{\sigma}\in W_{f,E/\mathbb{Q}}$ be such that
\[\bar {w}_{\sigma}|E=\sigma ,\quad\bar {w}_{\sigma\iota}=\bar {w}_{
\sigma}\tilde{\iota }.\]
Then (Langlands 1979; Milne and Shih 1982a), $\bar {b}(\tau )$ is defined by
\[\bar {b}(\tau )=\prod_{\sigma\in\Gal(E/\mathbb{Q})}\sigma\mu^E(\bar {
w}_{\sigma}\tilde{\tau}\bar {w}_{\sigma\tau}^{-1})\quad\mod S^E(E
).\]
Let $w_{\sigma}=\bar {w}_{\sigma^{-1}}^{-1}$; then $w_{\sigma}|E=
\sigma$ and $w_{\iota\sigma}=\tilde{\imath }w_{\sigma}$; moreover,
\[\bar {b}(\tau^{-1})^{-1}=\prod_{\sigma\in\Gal(E/\mathbb{Q})}\sigma
\mu^E(w_{\sigma^{-1}\tau}^{-1}\tilde{\tau }w_{\sigma^{-1}})=f^E(\tau
).\]
Thus, in the notation of Milne and Shih 1982a, 2.9, $f^E(\tau )=\bar{
\beta }(\tau )$.
\end{remark}

\section{Definition of $g^{E}$}

We complete the proof of (5.1) by showing that there is a universal $ g$
giving rise to all $g_{\Phi}$. For simplicity, we shall assume that $ E$ is
Galois over $\mathbb{Q}$ --- for a non Galois field, $g^{E}$ can be defined as
the norm of the element from the Galois closure.

\begin{proposition}
Let $E\subset \mathbb{C}$ be a CM-field.  There exists a unique map
$g^E:\Aut(\mathbb{C})\to S^E(\mathbb{A}_{f,E})/S^E(E)$ with the following property: for any CM-type
$(K,\Phi )$ whose reflex field is contained in $E$,
\[\rho_{\Phi}(g^E(\tau ))=g_{\Phi}(\tau )\]
in $T(\mathbb{A}_{f,E})/T(E)$, where $T=\Res_{K/\mathbb{Q}}(\mathbb{G}_
m)$.
\end{proposition}
\begin{proof}The uniqueness follows from (4.7).  For the existence, we need the
notion of a Hodge cycle.
For any variety $X$ over $\mathbb{C}$, write $H^s(X,\mathbb{Q})(r)=H^
s(X,(2\pi i)^r\mathbb{Q})$ (cohomology
with respect to the complex topology).  A {\em Hodge cycle on} $A$ is an element
$s\in H^{2p}(A^k,\mathbb{Q})(p)$, some $p,k$, that is of type $(p,p
)$, i.e.,  under the embedding
$(2\pi i)^p\mathbb{Q}\hookrightarrow \mathbb{C}$, $s$ maps into $H^{p
,p}\subset H^{2p}(X,\mathbb{C})$.  Recall that
$H^r(A^k,\mathbb{Q})=\bigwedge^r(\oplus^kH_1(A,\mathbb{Q})^{\vee})$, and so $\GL
(H_1(A,\mathbb{Q}))$ acts by transport of structure
on $H^r(A^k,\mathbb{Q})$.  The {\em Mumford-Tate group} $MT(A)$ of $
A$ is the largest
$\mathbb{Q}$-rational algebraic subgroup of $\GL(H_1(A,\mathbb{Q}))$ such that $
MT(A)(\mathbb{Q})$ is the set of
$\alpha\in\GL(H_1(A,\mathbb{Q}))$ for which there exists a $\nu (\alpha
)\in \mathbb{Q}^{\times}\text {\rm \ such that }\alpha s=\nu (\alpha
)^ps$
for any Hodge cycle $s$ on $A$ (of type $(p,p)$).
\begin{lemma}
Assume $A$ is of CM-type $(K,\Phi )$, where the reflex field of $
(K,\Phi )$ is contained
in $E$.  Then the image of $\rho_{\Phi}:S^E\to K^{\times}\subset\GL
(H_1(A,\mathbb{Q}))$ is equal to $MT(A)$.
\end{lemma}
\begin{proof}
Cf. Deligne 1982a, 3.4.
\end{proof}
Write $H^{2p}(A^k,\mathbb{A}_f)(p)=H^{2p}(A^k,\mathbb{Q})(p)\otimes \mathbb{A}_
f$.  Then there is a canonical
isomorphism
\[H^{2p}(A^k,\mathbb{A}_f)(p)\stackrel {\approx}{\to}\bigwedge^r(\oplus^
kV_f(A)^{\vee})\]
and so the action of $\Aut(\mathbb{C})$ on $V_f(A)$ gives rise to an action on
$H^{2p}(A^k,\mathbb{A}_f)(p)$.  We shall need to use the following important result fo
Deligne.
\begin{theorem}
Let $s\in H^{2p}(A^k,\mathbb{Q})(p)$ be a Hodge cycle on $A$, and let $
s_f$ be the image of $s$ in
$H^{2p}(A^k,\mathbb{A}_f)(p)$; then for any $\tau\in\Aut(\mathbb{C})$ there exists a Hodge cycle $
s_1$ on $\tau A$
whose image in $H^{2p}(A^k,\mathbb{A}_f)(p)$ is $\tau s_f$.
\end{theorem}
\begin{proof}
See Deligne 1982a.
\end{proof}
The cycle $s_1$ of the theorem is uniquely determined, and will be written $
\tau s$.
\begin{proposition}
With the notations of (8.2), there exists a $K$-linear isomorphism
$\alpha :H_1(A,E)\stackrel {\approx}{\to}H_1(\tau A,E)$ such that $
\alpha (s)=\nu (\alpha )^p\tau (s)$ for all Hodge cycles $s$ on $
A$
(of type $(p,p)$).
\end{proposition}
\begin{proof}
For any $\mathbb{Q}$-algebra $R$, let
\[P(R)=\{\alpha :H_1(A,R)\stackrel {\approx}{\to}H_1(\tau A,R)\mid
\alpha (s)=\nu (\alpha )^p\tau (s),\quad\text {\rm all }s\}.\]
Then $P(R)$ is either empty or is a principal homogeneous space over
$MT(A)(R)$.  Thus $P$ is either the empty scheme or is a principal
homogeneous space over $MT(A)$.  The existence of $\tau :H_1(A,\mathbb{A}_
f)\to H_1(\tau A,\mathbb{A}_f)$ in
$P(\mathbb{A}_f)$ shows that the latter is true.  It therefore corresponds to an
element of $H^1(\mathbb{Q},MT(A))$.  But $MT(A)_E\approx \mathbb{G}_m
\times\cdots\times \mathbb{G}_m$, and so $H^1(E,MT(A))=0$
by Hilbert's Theorem 90.
\end{proof}
Both (8.2) and (8.4) obviously also apply to products of abelian varieties of
CM-type.  Let $A=\prod A_{\Phi}$, where $\Phi$ runs through the CM-types on $
E$ and $A_{\Phi}$
is of type $(E,\Phi )$.  Then $\rho :S^E\stackrel {\approx}{\to}M
T(A)$.  Choose $\alpha$ as in (8.4).  Then
\[V_f(A)\otimes E\stackrel {\tau}{\to}V_f(\tau A)\otimes E\stackrel {
(\alpha\otimes 1)^{-1}}{\to}V_f(A)\otimes E\]
is an $\mathbb{A}_{f,E}$-linear isomorphism and sends a Hodge cycle $
s$ of type $(p,p)$ to
$\nu^ps$, some $\nu\in \mathbb{A}_{f,E}^{\times}$.  Therefore it is multiplication by an element
$g\in MT(A)(\mathbb{A}_{f,E})=S^E(\mathbb{A}_{f,E})$.  The class $g(\tau
)$ of $g$ in $S^E(\mathbb{A}_{f,E})/S^E(E)$ has
the properties required for (8.1).\end{proof}

The map $g:\Aut(\mathbb{C})\to S^{E}(\mathbb{A}_{f,E})/S^{E}(E)$ has the same
properties as those listed for $f$ in (7.3). In particular, $g(\tau)$ is fixed
by $\Gal(E/ \mathbb{Q})$. Set
\[
e(\tau)=\frac{g(\tau)}{f(\tau)}.
\]
Then $e(\tau)\in(S^{E}(\mathbb{A}_{f,E})/S^{E}(E))^{\Gal(E/\mathbb{Q})}$, and
it remains to show that it lies in $S^{E}(\mathbb{A}_{f})/S^{E}(\mathbb{Q})$
--- the next proposition completes the proof.

\begin{proposition}
$e(\tau )$ lies in $S^E(\mathbb{A}_f)/S^E(\mathbb{Q})$.
\end{proposition}
\begin{proof}
There is a cohomology sequence
\[0\to S^E(\mathbb{Q})\to S^E(\mathbb{A}_f)\to (S^E(\mathbb{A}_{f,E})/S^
E(E))^{\Gal(E/\mathbb{Q})}\to H^1(\mathbb{Q},S^E).\]
Thus, we have to show that the image $\gamma$ of $e(\tau )$ in $H^
1(\mathbb{Q},S^E)$ is zero.  But
$H^1(\mathbb{Q},S^E)\hookrightarrow\prod_{\ell ,\infty}H^1(\mathbb{Q}_{
\ell},S^E)$, as follows easily from (4.1), and the image of
$e(\tau )$ in $H^1(\mathbb{Q}_{\ell},S^E)$ is obviously zero for all finite $
\ell$.  It remains to check
that the image of $\gamma$ in $H^1(\mathbb{R},S^E)$ is zero.  Let
\[T=\{a\in\prod_{\text{\rm CM-types on }E}E^{\times}\mid a\cdot\iota
a\in \mathbb{Q}^{\times}\}\qquad(\text{ torus over }\mathbb{Q}).\]
\begin{lemma}
The image of $\gamma$ in $H^1(\mathbb{Q},T)$ is zero.
\end{lemma}
\begin{proof}
In the proof of (3.6) it shown that the image of $e$ in $T(\mathbb{A}_{
f,E})/T(E)$ lifts to
an element $\varepsilon\in T(\mathbb{A}_f)$.  The image of $\gamma$ in $
H^1(\mathbb{Q},T)$ is represented by the
cocycle $\sigma\mapsto\sigma\varepsilon -\varepsilon =0$.
\begin{lemma}
The map $H^1(\mathbb{R},S^E)\to H^1(\mathbb{R},T)$ is injective.
\end{lemma}
\begin{proof}
There is a norm map $a\mapsto a\cdot\iota a:T\twoheadrightarrow \mathbb{G}_
m$, and we define $ST$ and
$SMT(A)$ to make the rows in
\[
\begin{diagram}
1&\rTo&SMT(A)&\rTo&MT(A)&\rTo&{\mathbb{G}}_m&\rTo&1\\
&&\dInto&&\dInto&&\|\\
1&\rTo&ST&\rTo&T&\rTo&{\mathbb{G}}_m&\rTo&1
\end{diagram}
\]
exact.  (Here $A=\prod A_{\Phi}$.)  This diagram gives rise to an exact commutative
diagram
\[
\begin{diagram}
{\mathbb{R}}^{\times}&\rTo&H^1({\mathbb{R}},SMT)&\rTo&H^1({\mathbb{R}},MT)&\rTo&0\\
\|&&\dTo&&\dTo\\
{\mathbb{R}}^{\times}&\rTo&H^1({\mathbb{R}},ST)&\rTo&H^1({\mathbb{R}},T)&\rTo&0.
\end{diagram}
\]
Note that $ST$ (and hence $SMT$) are anisotropic over $\mathbb{R}$; hence,
$H^1(\mathbb{R},SMT)=SMT(\mathbb{C})_2$  and $H^1(\mathbb{R},ST)=ST(\mathbb{C}
)_2$, and so
$H^1(\mathbb{R},SMT)\hookrightarrow H^1(\mathbb{R},SMT)$.  The five-lemma now completes the proof.
\end{proof}\end{proof}\end{proof}

See also Milne and Shih, 1982b, \S 5.

\begin{remark}
It seems to be essential to make use of Hodge cycles, and consequently
Shimura varieties (which are used in the proof of (8.3)), in order to show
the $e_{\Phi}(\tau )$ have the correct functorial properties.  Note that Shih (1976) also
needed to use Shimura varieties to prove his case of the theorem.
\end{remark}

\section{Re-statement of the Theorem}

The following statement of the main theorem of complex multiplication first
appeared (as a conjecture) in Milne and Shih 1979.

\begin{theorem}
Let $A$ be an abelian variety of CM-type $(K,\Phi )$; let $\tau
\in\Aut(\mathbb{C})$, and let
$f\in f(\tau )$.  Then
\begin{enumerate}
\item\ $\tau A$ is of type $(K,\tau\Phi )$;
\item\ there is an $K$-linear isomorphism $\alpha :H_1(A,E)\to H_
1(\tau A,E)$ where $E$ is the
reflex field of $(K,\Phi )$, such that
\begin{enumerate}
\item\ $\alpha (s)=\nu (\alpha )^p\tau (s)$, for all Hodge cycles $
s$ on $A$, where $\nu (\alpha )\in \mathbb{Q}^{\times}$ and $2p$
is the degree of $s$;
\item\
\[\begin{diagram}[heads=LaTeX,size=3em]
V_f(A)\otimes{}E&\rTo{\rho_\Phi(f)}&V_f(A)\otimes{}E\\
&\rdTo{}{\tau}&\dTo{\alpha\otimes{1}}\\
&&V_f(\tau{A})\otimes{}E\end{diagram}\]
commutes (note that $\rho_{\Phi}(f)\in \mathbb{A}_{f,K\otimes E}^{\times}
)$.
\end{enumerate}
\end{enumerate}
\end{theorem}
\begin{proof}
The theorem is true (by definition) if $f(\tau )$ is replaced by $
g(\tau )$, but we have
shown that $g(\tau )=f(\tau )$.
\end{proof}

\begin{remark}
Let $T$ be a torus such that
\[MT(A)\subset T\subset \{a\in K^{\times}\mid a\cdot\iota a\in \mathbb{Q}^{
\times}\}\]
and let $h$ be the homomorphism defining the Hodge structure on $
H_1(A,\mathbb{R})$.
Then the Shimura variety $Sh(T,\{h\})$ is, in a natural way, a moduli scheme,
and the (new) main theorem of complex multiplication gives a description of
the action of $\Aut(\mathbb{C})$ on $Sh(T,\{h\})$ (see Milne and Shih 1982, \S 6).
\end{remark}

\begin{remark}
Out of his study of the zeta functions of Shimura varieties, Langlands (1979)
was led to a conjecture concerning the conjugates of Shimura varieties.
The conjecture is trivial for the Shimura varieties associated with tori, but
in Milne and Shih 1982b it is shown that for groups of symplectic
similitudes the conjecture is equivalent to (9.1).  It is also shown (ibid.)
that the validity of the conjecture for a Shimura variety $Sh(G,X
)$ depends
only on $(G^{\text {\rm der}},X^{+})$.  Thus (ibid.)  similar methods to those used in Deligne
1979a can be used to prove Langlands's conjecture for exactly those Shimura
varieties for which Deligne proves the existence of canonical models
in that article.
\end{remark}

\section{The Taniyama Group}

By an extension of $\Gal(\mathbb{Q}^{\al}/\mathbb{Q})$ by $S^{E}$ with
finite-ad\`elic splitting, we mean an exact sequence
\[
1\to S^{E}\to T^{E}\overset{\pi^{E}}{\to}\Gal(\mathbb{Q}^{\al}/\mathbb{Q}
)\to1
\]
of pro-algebraic groups over $\mathbb{Q}$ ($\Gal(\mathbb{Q}^{\al}/\mathbb{Q}
)$ is to be regarded as a constant pro-algebraic group) together with a
continuous homomorphism $sp^{E}:\Gal(\mathbb{Q}^{\al}/\mathbb{Q})\to
T^{E}(\mathbb{A}_{f})$ such that $ sp^{E}\circ\pi^{E}=\id$. We always assume
that the action of $\Gal(\mathbb{Q}^{\al}/\mathbb{Q})$ on $S^{E}$ given by the
extension is the natural action. Assume $E\subset\mathbb{C}$ is Galois over
$\mathbb{Q}$, and a CM-field.

\begin{proposition}
\begin{enumerate}
\item\ Let $(T^E,sp^E)$ be an extension of $\Gal(\mathbb{Q}^{\al}/\mathbb{Q}
)$ by $S^E$ with finite-ad\`elic
splitting.  Choose a section $a^E:\Gal(\mathbb{Q}^{\al}/\mathbb{Q})\to
(T^E)_E$ that is a morphism of
pro-algebraic groups.  Define $h(\tau )\in S^E(\mathbb{A}_{f,E})/S^
E(E)$ to be the class of
$sp^E(\tau )\cdot a^E(\tau )^{-1}$.
\begin{enumerate}
\item\ $h(\tau )$ is well-defined;
\item\ $\sigma h(\tau )=h(\tau )$, $\sigma\in\Gal(E/\mathbb{Q})$;
\item\ $h(\tau_1\tau_2)=h(\tau_1)\cdot\tilde{\tau}_1h(\tau_2)$, $
\tau_1,\tau_2\in\Gal(\mathbb{Q}^{\al}/\mathbb{Q})$;
\item\ $h$ lifts to a continuous map $h':\Gal(\mathbb{Q}^{\al}/\mathbb{Q}
)\to S^E(\mathbb{A}_{f,E})$ such that the
map $(\tau_1,\tau_2)\mapsto d_{\tau_1,\tau_2}\stackrel {\text{\rm df}}{
=}h'(\tau_1)\cdot\tilde{\tau}_1h'(\tau_2)\cdot h'(\tau_1\tau_2)^{
-1}$ is locally constant.
\end{enumerate}
\item\ Let $h:\Gal(\mathbb{Q}^{\al}/\mathbb{Q})\to S^E(\mathbb{A}_{f,E}
)/S^E(E)$ be a map satisfying conditions (i),
(ii), (iii), (iv); then $h$ arises from a unique extension of $\Gal
(\mathbb{Q}^{\al}/\mathbb{Q})$ by $S^E$
with finite-ad\`elic splitting.
\end{enumerate}
\end{proposition}
\begin{proof}
Easy; see Milne and Shih 1982a, \S 2.
\end{proof}

Let $S=\plim S^{E}$, where $E$ runs through the CM-fields $E\subset\mathbb{C}$
that are Galois over $\mathbb{Q}$. By an extension of $\Gal(\mathbb{Q}%
^{\al}/\mathbb{Q} )$ by $S$ with finite-ad\`elic splitting, we mean a
projective system of extensions of $\Gal(\mathbb{Q}^{\al} /\mathbb{Q})$ by
$S^{E}$ with finite-ad\'elic splitting, i.e., a family
\[
\begin{CD}
1@>>>S^{E_1}@>>>T^{E_1}@>>>\Gal({\mathbb{Q}}^{\text{al}}/{\mathbb{Q}})@>>>1\\
@.@VV{\Nm_{E_1/E_2}}V@VV{\Nm_{E_1/E_2}}V@VV{\id}V\\
1@>>>S^{E_2}@>>>T^{E_2}@>>>\Gal({\mathbb{Q}}^{\text{al}}/{\mathbb{Q}})@>>>1
\end{CD}
\]

\[
\begin{diagram}[heads=LaTeX]
T^{E_1}({\mathbb{A}}_f)\\
&\luTo{sp^{E_1}}\\
\dTo{\Nm_{E_1/E_2}}&&\Gal({\mathbb{Q}}^{\text{al}}/{\mathbb{Q}})\\
&\ldTo{}{sp^{E_2}}\\
T^{E_2}({\mathbb{A}}_f)
\end{diagram}
\]
of commutative diagrams.

\begin{theorem}
Let $T_1$ and $T_2$ be two extensions of $\Gal(\mathbb{Q}^{\al}/\mathbb{Q}
)$ by $S$ with finite-ad\`elic
splittings.  Assume:
\begin{enumerate}
\item\ for each $E$, and $i=1,2$, there exists a commutative diagram
\[\begin{diagram}
1&\rTo&S^E&\rTo&{}_{E}T_i^E&\rTo&\Gal({\mathbb{Q}}^{\text{al}}/{\mathbb{Q}})&\rTo&1\\
&&\|&&\dTo&&\dOnto\\
1&\rTo&S^E&\rTo&M^E&\rTo&\Gal({\mathbb{Q}}^{\text{al}}/E)^{\text{ab}}&\rTo&1\\
\end{diagram}\]
compatible with the finite-ad\`elic splittings, where $_ET^E$ is the inverse image
of $\Gal(\mathbb{Q}^{\al}/\mathbb{Q})$ in $T^E$ and the lower row is the extension constructed by
Serre (1968, II).
\item\ for each $\tau\in\Gal(\mathbb{Q}^{\al}/\mathbb{Q})$, $\pi_1^{-
1}(\tau )\approx\pi_2^{-1}(\tau )$ as principal homogeneous
spaces over $S^E$.
\item\ $sp^E(\iota )\in T_i^E(\mathbb{Q})$, $i=1,2$.
\end{enumerate}
Then there is a unique family of isomorphisms $\phi^E:T_1^E\to T_
2^E$ making the
following diagrams commute:
\[\begin{CD}
1@>>>S^E@>>>T_1^E@>>>\Gal({\mathbb{Q}}^{\text{al}}/{\mathbb{Q}})@>>>1\\
@.@VV{\id}V@VV{\phi^E}V@VV{\id}V\\
1@>>>S^E@>>>T_2^E@>>>\Gal({\mathbb{Q}}^{\text{al}}/{\mathbb{Q}})@>>>1
\end{CD}\]
\[\begin{CD}
T_1^{E_1}@>{\Nm_{{E_1}/{E_2}}}>>T_1^{E_2}\\
@VV{\phi^{E_1}}V@VV{\phi^{E_2}}V\\
T_2^{E_2}@>{\Nm_{{E_1}/{E_2}}}>>T_2^{E_2}
\end{CD}\]
\[\begin{CD}
T_1^E({\mathbb{A}}_f)@<{sp_1^E}<<\Gal({\mathbb{Q}}^{\text{al}}/{\mathbb{Q}})\\
@VVV@VV{\id}V\\
T_2^E({\mathbb{A}}_f)@<{sp_2^E}<<\Gal({\mathbb{Q}}^{\text{al}}/{\mathbb{Q}})
\end{CD}\]
\end{theorem}
\begin{proof}
Let $(h_1^E)$ and $(h_2^E)$ be the families of maps corresponding as in (10.1a) to $
T_1$
and $T_2$.  The hypotheses of the theorem imply that the family $
(e^E)$, where
$e^E=h_1^E/h_2^E$, satisfies the hypotheses of (6.1).  Thus $h_1^
E=h_2^E$ for all $E$, and
we apply (10.1b).
\end{proof}

\begin{definition}
The extension corresponding to the family of maps $(f^E)$ (rather,
$\tau\mapsto f^E(\tau^{-1})^{-1}$)  defined in (7.3) is called the Taniyama group.
\end{definition}

\begin{remark}
In (1982b), Deligne proves the following:
\begin{enumerate}
\item\ let $T'$ be the group associated with the Tannakian category of motives
over $\mathbb{Q}$ generated by Artin motives and abelian varieties of potential
CM-type; then $T'$ is an extension of $\Gal(\mathbb{Q}^{\al}/\mathbb{Q}
)$ by $S$ with finite-ad\`elic
splitting in the sense defined above.  (From a more naive point of view, $
T'$
is the extension defined by the maps $(g^E)$ of $\S 8$.)
\item\ Theorem 10.2, by essentially the same argument as we have given in
\S 6, except expressed directly in terms of the extensions rather than
cocycles.
\end{enumerate}
These two results combine to show that the motivic Galois group is
isomorphic to the explicitly constructed Taniyama group (as extensions with
finite-ad\`elic splitting).
This can be regarded as another statement of the (new) main theorem of
complex multiplication.
Note however that without the Taniyama group, Deligne's result says
very\footnote{It is only a uniqueness result:  it says that there is at most
extension consistent with the theorem of Shimura and Taniyama; Langlands
wrote down an {\em explicit\/} extension with this
property.} little.  This is why I have included Langlands as one of the main
contributors\footnote{Probably Shih and Tate should also be included.} to the proof of (1.1) even though he never explicitly considered
abelian varieties with complex multiplication (and neither he nor Deligne
explicitly considered a statement like (1.1)).
\end{remark}

\section{Zeta Functions}

\begin{lemma}
There exists a commutative diagram
\[
\begin{diagram}
T^E({\mathbb{C}})&\lTo{sp_{\infty}^E}&W_{\mathbb{Q}}\\
\uInto&&\dTo\\
T^E({\mathbb{Q}})&\rTo&\Gal({\mathbb{Q}}^{\text{al}}/{\mathbb{Q}})
\end{diagram}
\]
where $W_{\mathbb{Q}}$ is the Weil group of $\mathbb{Q}$ and $T$ is the Taniyama group.
\end{lemma}
\begin{proof}
Easy; see Milne and Shih 1982a, 3.17.
\end{proof}

\begin{theorem}
Let $A$ be an abelian variety over $\mathbb{Q}$ of potential CM-type $
(K,\Phi )$.  Let $E$ be a
CM-field containing the reflex field of $(K,\Phi )$.  Then there exists a
representation $\rho :T^E\to\Aut(H_1(A_{\mathbb{C}},\mathbb{Q}))$ such that
\begin{enumerate}
\item\ $\rho_f\stackrel {\text{\rm df}}{=}\rho\circ sp^E:\Gal(\mathbb{Q}^{\al}
/\mathbb{Q})\to\Aut(V_f(A))$ describes the action of
$\Gal(\mathbb{Q}^{\al}/\mathbb{Q})$ on $V_f(A)$;
\item\ $L(s,A/\mathbb{Q})=L(s,\rho_{\infty})$ where $\rho_{\infty}=
\rho\circ sp_{\infty}^E$ is a complex representation of
$W_E$.
\end{enumerate}
\end{theorem}
\begin{proof}
The existence of $\rho$ is obvious from the interpretation of $T$ as the motivic
Galois group $M$ (see 10.4a) or, more naively, as the extension corresponding
to $(g^E(\tau^{-1})^{-1})$.
\end{proof}

\begin{remark}
The proof of (11.2) does not require the full strength of Deligne's results,
and in fact is proved by Deligne (1979b).  Subsequently Yoshida (1981) found
another proof that $L(s,A/\mathbb{Q})=L(s,\rho_{\infty})$ for {\em some\/} complex representation $
\rho_{\infty}$
of $W_{\mathbb{Q}}$.  When $\Gal(\mathbb{Q}^{\al}/\mathbb{Q})$ stabilizes $
K\subset\End(A_{\mathbb{Q}^{\al}})\otimes \mathbb{Q}$, this last result
was proved independently by Milne (1972) (all primes) and Shimura (1971) (good primes only).
\end{remark}

\centerline{\textbf{References}}

Deligne, P., Vari\'et\'es de Shimura: interpretation modulaire, et techniques
de construction de mod\`eles canoniques. Proc. Symp. Pure Math., A.M.S., 33
(1979a), part 2, 247--290.

Deligne, P., Letter to Langlands, 10 April, 1979b.

Deligne, P., Hodge cycles on abelian varieties, 1982a. In: Deligne et al.1982,
pp 9--100.

Deligne, P., Motifs et groupe de Taniyama, 1982b, (completion of Deligne
1979b). In: Deligne et al. 1992, pp 261--279.

Deligne, P., Milne, J.S., Ogus, A., Shih, Kuang-yen, Hodge Cycles, Motives,
and Shimura Varieties, Lecture Notes in Mathematics 900, Springer, 1982.

Langlands, R., Automorphic representations, Shimura varieties, and motives.
Ein M\"archen. Proc. Symp. Pure Math., A.M.S., 33 (1979), part 2, 205--246.

Milne, J.S., On the arithmetic of abelian varieties, Invent. Math. 17 (1972), 177-190.

Milne, J.S., and Shih, Kuang-yen, Shimura varieties: Conjugates and the action
of complex conjugation, 154pp, 1979. (Published in revised form as Milne and
Shih 1981, 1982a, 1982b.)

Milne, J.S., and Shih, Kuang-yen, Langlands's construction of the Taniyama
group, 1982a. In: Deligne et al 1982, pp 229--260.

Milne, J.S., and Shih, Kuang-yen, Conjugates of Shimura varieties, 1982b. In:
Deligne et al 1982, pp 280--356.

Serre, J-P., Sur les groupes de congruence des vari\'et\'es ab\'eliennes, Izv.
Akad. Nauk. SSSR 28 (1964), 3--18.

Serre, J-P., Abelian $l$-adic Representations and Elliptic Curves, Benjamin, 1968.

Shih, Kuang-yen, Anti-holomorphic automorphisms of arithmetic automorphic
function fields, Annals of Math. 103 (1976), 81--102.

Shimura, G., Arithmetic Theory of Automorphic Functions, Princeton U. P., 1971.

Shimura, G., On the zeta-function of an abelian variety with complex
multiplication, Annals of Math. 94 (1971), 504--533.

Tate, J., Number theoretic background, in \textit{Automorphic forms,
representations and $L$-functions (Proc. Sympos.Pure Math., Oregon State
Univ., Corvallis, Ore., 1977), Part 2}, 3--26, Proc. Sympos. Pure Math.,
XXXIII, Amer. Math. Soc., Providence, R.I., 1979.

Tate, J., On conjugation of abelian varieties of CM-type, 8pp, April 1981.

Yoshida, H., Abelian varieties with complex multiplication and representations
of the Weil groups. Annals of Math. 114 (1981), 87--102.

September 19, 1981.

\newpage\centerline{\textbf{\Large Addendum (June 1998)}}

The sections of the Addendum are largely independent.

\section[Origins]{The Origins of the Theory of Complex Multiplication for Abelian
Varieties of Dimension Greater Than One}

On this topic, one can not do better than to quote Weil's commentary
(\OE uvres, Vol II, pp 541--542) on his articles in the
Proceedings\footnote{This is the same conference where Taniyama gave his
somewhat enigmatic statement of the Taniyama conjecture.} of the International
Symposium on Algebraic Number Theory, held in Tokyo and Nikko, September
8--13, 1955.

\begin{quotation}
Comme contribution au colloque, j'apportais quelques id\'ees que je croyais
neuves sur l'extension aux vari\'et\'es ab\'eliennes de la th\'eorie classique
de la multiplication complexe. Comme chacun sait, Hecke avait eu l'audace,
stup\'efiante pour l'\'epoque, de s'attaquer \`a ce probl\`eme d\`es 1912; il
en avait tir\'e sa th\`ese, puis avait pouss\'e son travail assez loin pour
d\'ecouvrir des ph\'enom\`enes que lui avaient paru inexplicables, apr\`es
quoi il avait abandonn\'e ce terrain de recherche dont assur\'ement
l'exploration \'etait pr\'ematur\'ee. En 1955, \`a la lumi\`ere des progr\`es
effectu\'es en g\'eom\'etrie alg\'ebrique, on pouvait esp\'erer que la
question \'etait m\^ure.

Elle l'\'etait en effet; \`a peine arriv\'e \`a Tokyo, j'appris que deux
jeunes japonais venaient d'accomplir sur ce m\^eme sujet des progr\`es
d\'ecisifs. Mon plaisir \`a cette nouvelle ne fut un peu temp\'er\'e que par
ma crainte de n'avoir plus rien \`a dire au colloque. Mais il apparut
bient\^ot, d'abord que Shimura et Taniyama avaient travaill\'e
ind\'ependamment de moi et m\^eme ind\'ependamment l'un de l'autre, et surtout
que nos r\'esultats \`a tous trois, tout en ayant de larges parties communes,
se compl\'etaient mutuellement. Shimura avait rendu possible la r\'eduction
modulo $\mathfrak{p}$ au moyen de sa th\'eorie des intersections dans les
vari\'et\'es d\'efinies sur un anneau local (\emph{Am. J. of Math.\/} 77
(1955), pp. 134--176); il s'en \'etait servi pour l'\'etude de vari\'et\'es
ab\'eliennes \`a multiplication complexe, bien qu'initialement, \`a ce qu'il
me dit, il e\^ut plut\^ot eu en vue d'autres applications. Taniyama, de son
c\^ot\'e, avait concentr\'e son attention sur les fonctions z\^eta des
vari\'et\'es en question et principalement des jacobiennes, et avait
g\'en\'eralis\'e \`a celles-ci une bonne partie des r\'esultats de Deuring sur
le cas elliptique. Quant \`a ma contribution, elle tenait surtout \`a l'emploi
de la notion de ``vari\'et\'e polaris\'ee''; j'avais choisi ce terme, par
analogie avec ``vari\'et\'es orient\'ees'' des topologues, pour d\'esigner une
structure suppl\'ementaire qu'on peut mettre sur une vari\'et\'e compl\`ete et
normale quand elle admet un plongement projectif. Faute de cette structure, la
notion de modules perd son sens.

Il fut convenu entre nous trois que je ferais au colloque un expos\'e
g\'en\'eral ([Weil 1956b]) \'esquissant \`a grands traits l'ensemble des
r\'esultats obtenus, expos\'e qui servirait en m\^eme temps d'introduction aux
communications de Shimura et de Taniyama; il fut entendu aussi que par la
suite ceux-ci r\'edigeraient le tout avec des d\'emonstrations d\'etaill\'ees.
Leur livre a paru en 1961 sous le titre \emph{Complex multiplication of
abelian varieties and its } \emph{application to number theory\/} (Math. Soc.
of Japan, Tokyo); mais Taniyama \'etait mort tragiquement en 1958, et Shimura
avait d\^u l'achever seul.

D'autre part, tout en restant loin des r\'esultats de Taniyama sur les
fonctions z\^eta des vari\'et\'es ``de type CM'' (comme on dit \`a pr\'esent),
j'avais aper\d cu le r\^ole que devaient jouer dans cette th\'eorie certains
caract\`eres de type $(A_{0})$'', ainsi que les caract\`eres \`a values
$\mathfrak{P}$-adiques qu'ils permettent de d\'efinir (cf. [Weil 1955b], p6).
Je trouvai l\`a une premi\`ere explication du ph\'enom\`ene qui avait le plus
\'etonn\'e Hecke; il consiste en ce que, d\`es la dimension $2$, les modules
et les points de division des vari\'et\'es de type CM d\'efinissent en
g\'eneral des extensions ab\'eliennes, non sur le corps de la multiplication
complexe, mais sur un autre qui lui est associ\'e. Ce sujet a \'et\'e repris
et plus amplement d\'evelopp\'e par Taniyama (\emph{J. Math. Soc. }
\emph{Jap.\/} 9 (1957), pp. 330--366); cf. aussi [Weil 1959].
\end{quotation}

As mentioned in the text, the first theorem extending the Main Theorem of
Complex Multiplication to automorphisms not fixing the reflex field was that
of Shih 1976. This theorem of Shih was used in Milne and Shih 1981 to give an
explicit description of the involution defined by complex conjugation on the
points of Shimura variety whose reflex field is real (Conjecture of Langlands
1979, p 234).

Apparently, it was known to Grothendieck, Serre, and Deligne in the 1960s that
the conjectural theory of motives had as an explicit consequence the existence
of a Taniyama group --- these ideas inform the presentation in Chapters I and
II of Serre 1968 --- but they were unable to construct such a group. It was
not until 1977, when Langlands's efforts to understand the conjugates of
Shimura led him to define his cocycles, that the group could be constructed
(and it was Deligne who recognized that Langlands's cocycles answered the
earlier problem).

The rest of the story is described in the text of the article.

\section{Zeta Functions of Abelian Varieties of CM-type}

In this section I explain the elementary approach (Milne 1972), not using the
theorems in the first part of this article, to the zeta function of abelian
varieties of CM-type.

First some terminology: For abelian varieties $A$ and $B$ over a field $ k$,
$\Hom(A,B)$ denotes the group of homomorphisms $A\to B$ defined over $ k$, and
$\Hom^{0}(A,B)=\Hom(A,B)\otimes_{\mathbb{Z}}\mathbb{Q}$. Similar notations are
used for endomorphisms. An abelian variety over $k$ is \emph{simple\/} if it
contains no nonzero proper abelian subvariety defined over $k$, and it is
\emph{absolutely simple\/} if it is simple\footnote{An older terminology,
based on Weil's Foundations, uses ``simple'' where we use ``absolutely
simple'', see, for example, Lang 1983 or Shimura 1998.} over $k^{\al}$. An
abelian variety $ A$ over a field $k$ of characteristic zero is said to be of
\emph{CM-type\/} if its Mumford-Tate group is a torus. Thus, $A$ is of CM-type
if, for each simple isogeny factor $B$ of $A_{k^{\al}}$, $\End^{0}(B)$ is a
CM-field of degree $2\dim B$ over $\mathbb{Q}$. For an abelian variety $A$
over a number field $k\subset\mathbb{C}$ and finite prime $ v$ of $k$, the
polynomial
\[
P_{v}(A,T)=\det(1-F_{v}T|V_{\ell}(A)^{I_{v}})
\]
where $\ell$ is any prime number different from the characteristic of the
residue field at $v$, $I_{v}$ is the inertia group at a prime $v^{\prime}|v$,
and $ F_{v}$ is a Frobenius element in the quotient of the decomposition group
at $ v^{\prime}$ by $I_{v}$ --- it is known that $P_{v}(A,T)$ is independent
of the choice of $ \ell$, $v^{\prime}$, and $F_{v}$. Finally, the zeta
function of $A$ is
\[
\zeta(A,s)=\prod_{v}\frac1{P_{v}(\mathbb{N}(v)^{-s})}%
\]
where $v$ runs over \emph{all\/} finite primes of $k$ and $\mathbb{N} (v)$ is
the order of the residue field at $v$. Clearly $\zeta(A,s)$ depends only on
the isogeny class of $ A$, and if $A$ is isogenous to $A_{1}\times\cdots\times
A_{m}$, then $\zeta(A,s)=\prod_{ i=1}^{m}\zeta(A_{i},s)$.

\subsection{Case that all endomorphisms of $A$ are defined over $ k$}

In this subsection, $A$ is an abelian variety of CM-type such that
$\End(A)=\End(A_{k^{\al}})$. Because $\zeta(A,s)$ depends only on the isogeny
class of $ A$, we may suppose that $A$ is isotypic, i.e., that it is isogenous
to a power of a simple abelian variety. Then there exists a CM-field
$K\subset\End^{0}(A)$ of degree $2\dim A$ over $\mathbb{Q}$.

The tangent space $T$ to $A$ is a finite-dimensional vector space over $
\mathbb{Q}$ on which both $k$ and $K$ act. Since $K$ acts $k$-linearly, the
actions commute. An element $\alpha\in k^{\times}$ defines an automorphism of
$T$ viewed as $ K$-vector space, whose determinant we denote $\psi_{0}%
(\alpha)$. Then $\psi_{0}:k^{\times} \to K^{\times}$ is a homomorphism. Let
$\mathbb{I}_{k}$ denote the group of id\`eles of $k$.

\begin{theorem}
There exists a unique homomorphism
\[\varepsilon :\mathbb{I}_k\to K^{\times}\]
such that
\begin{enumerate}
\item\ the restriction of $\varepsilon$ to $k^{\times}$ is $\psi_
0$;
\item\ the homomorphism $\varepsilon$ is continuous, in the sense that its kernel is
open in $\mathbb{I}_k$;
\item\ there is a finite set $S$ of primes of $k$, including those where $
A$ has
bad reduction, such that for all finite primes $v\notin S$, $\varepsilon$ maps any prime
element at $v$ to $F_v$.
\end{enumerate}
\end{theorem}
\begin{proof}
This is a restatement of the Theorem of Shimura and Taniyama (1961, p148)
--- see Serre and Tate 1968, Theorem 10.
\end{proof}

There is a unique continuous homomorphism $\chi:\mathbb{I}_{k}\to(K
\otimes_{\mathbb{Q}}\mathbb{R})^{\times}$ that is trivial on $k^{\times}$ and
coincides with $\varepsilon$ on the group $\mathbb{I}_{k}^{\infty}$ of
id\`eles whose infinite component is $1$ (ib. p513). For each $\sigma
:K\to\mathbb{C}$, let $ \chi_{\sigma}$ be the composite
\[
\mathbb{I}_{k}\overset{\chi}{\to}(K\otimes_{\mathbb{Q}}\mathbb{R})^{\times
}\overset{ \sigma\otimes1}{\to}\mathbb{C}^{\times}.
\]
It is continuous and trivial on $k^{\times}$, that is, it is a Hecke character
in the broad sense (taking values in $\mathbb{C}^{\times}$ rather than the
unit circle).

\begin{theorem}
The zeta function of $A$,
\[\zeta (A,s)=\prod_{\sigma :K\hookrightarrow \mathbb{C}}L(s,\chi_{
\sigma}).\]
\end{theorem}
\begin{proof}
This is proved in Shimura and Taniyama 1961 except for the factors
corresponding to a finite set of primes, and for all primes in Serre and
Tate 1968.
\end{proof}

\subsection{General Case}

We now explain how to extend these results to abelian varieties that are of
CM-type, but whose endomorphisms are not defined over the given field of definition.

Let $k$ be a field of characteristic zero, and let $A$ be an abelian variety
over a finite extension $k^{\prime}$ of $k$. The restriction of scalars
$\Res_{ k^{\prime}/k}A$ of $A$ to $k$ is the variety $A_{*}$ over $k$
representing the functor of $k$-algebras, $ R\mapsto A(R\otimes_{k}k^{\prime
})$. For any finite Galois extension $\bar{k}$ of $k$ containing $k^{\prime}$,
there is a canonical isomorphism
\[
P:A_{*\bar{k}}\overset{\approx}{\to}\prod_{\sigma\in\Hom_{k}(k^{\prime}
,\bar{k})}\sigma A.
\]

\begin{lemma}
Let $k$ be a number field, and let $A_{*}$ be the abelian variety over $
k$ obtained
by restriction of scalars from an abelian variety $A$ over a finite extension
$k'$ of $k$.  Then $\zeta (A_{*},s)=\zeta (A,s)$.
\end{lemma}
\begin{proof}
It is immediate from the definition of $A_{*}$ that $V_{\ell}(A_{
*})$ is the
$\Gal(k^{\al}/k)$-module induced from the $\Gal(k^{\al}/k')$-module $
V_{\ell}(A)$.  This implies
the statement.  (See Milne 1972, Proposition 3.)
\end{proof}

\begin{lemma}
Let $A$ be an abelian variety over a field $k$, and let $k'$ be a finite Galois
extension of $k$ of degree $m$ and Galois group $G$.  Suppose that there exists a
$\mathbb{Q}$-subalgebra $R\subset\End^0(A_{k'})$ such that $R^G$ is a field and $
[R:R^G]=m$.  Then
$\Res_{k'/k}A_{k'}$ is isogenous to $A^m$.
\end{lemma}
\begin{proof}
Let $\alpha_1,\ldots ,\alpha_m$ be an $R^G$-basis for $R$ over $R^
G$, and let $\phi :A_{k'}^m\to A_{k'}^m$ be the
homomorphism $(\sigma_i\alpha_j)_{1\leq i,j\leq m}$, where $G=\{\sigma_
1,\ldots ,\sigma_m\}$.  Then $\phi$ is an isogeny.  When
we identify the second copy of $A_{k'}^m$ with $\prod\sigma_iA_{k'}$ and compose $
\phi$ with $P^{-1}$,
we obtain an isogeny $A_{k'}^m\to A_{*k'}$ that is invariant under $
G$, and hence defined
over $k$ (ibid.  Theorem 3).
\end{proof}

\begin{example}
Let $A$ be a simple abelian variety over a field $k$.  Let $R$ be the centre of
$\End^0(A_{k^{\al}})$, and let $k'$ be the smallest field containing $
k$ and such that all
elements of $R$ are over defined over $k'$.  Then $A$, $k'$, and $
R$ satisfy the
hypotheses of Lemma 13.4 (ibid.  p186), and so $A^m$ is isogenous to $
(A_{k'})_{*}$.
Hence, when $k$ is a number field
\[\zeta (A,s)^m=\zeta (A_{k'},s).\]
\end{example}

\begin{example}
Let $A$ be an abelian variety over a number field $k$ that, over $
\mathbb{C}$, becomes of
CM-type $(K,\Phi )$ for some field $K$.  Assume that $K$ is stable under the action
of $\Gal(k^{\al}/k)$ on $\End^0(A_{k^{\al}})$, and let $k'$ be the smallest field containing $
k$
such that all elements of $K$ are defined over $k'$.  Then $A$, $
k'$, and $K$ satisfy
the hypotheses of Lemma 13.4,  and so
\[\zeta (A,s)^m=\zeta (A_{k'},s)\stackrel {(13.2)}{=}\prod_{\sigma
\in\Sigma}L(s,\chi_{\sigma}),\quad\Sigma =\Hom(K,\mathbb{C}).\]
In this case, we can improve the result.  The group $G=\Gal(k'/k)$ acts
faithfully on $K$, and a direct calculation shows that $L(s,\chi_{
\sigma\circ\tau})=L(s,\chi_{\sigma})$ for
all $\sigma\in\Sigma$ and $\tau\in G$.  Therefore,
\[\prod_{\sigma\in\Sigma}L(s,\chi_{\sigma})=(\prod_{\sigma\in
\Sigma /G}L(s,\chi_{\sigma}))^m.\]
We can take an $m$th root, and obtain
\[\zeta (A,s)=\prod_{\sigma\in\Sigma /G}L(s,\chi_{\sigma}).\]
\end{example}

Now we consider the general case. Let $A$ be an abelian variety of CM-type
over a number field $k$. As noted earlier, we may suppose $A$ to be simple.
Then $\zeta(A,s)^{m}=\zeta(A_{k^{\prime}},s)$ where $k^{\prime}$ is the
smallest field containing $ k$ over which all endomorphisms in the centre of
$\End^{0}(A_{k^{\al}})$ are defined. Replacing $A/k$ with $A_{k^{\prime}%
}/k^{\prime}$, we may suppose that the endomorphisms in the centre of
$\End^{0}(A_{k^{\al}})$ are defined over $k$, and we may again suppose that $
A$ is simple. Then $A_{k^{\al}}$ is isotypic, and, for example, if it is
simple, we can apply (13.5) to obtain the zeta function of $A$.

\section{Hilbert's Twelfth Problem}

This asks for

\begin{quotation}
\ldots\ those functions that play for an arbitrary algebraic number field the
role that the exponential function plays for the field of rational numbers and
the elliptic modular functions play for an imaginary quadratic number field.
\end{quotation}

The classical result, referred to by Hilbert, can be stated as follows: for
any quadratic imaginary field $E$, the maximal abelian extension of $ E$ is
obtained by adjoining to it the moduli of elliptic curves and their torsion
points with complex multiplication by $E$.

As Weil observed (see \S 12), in dimension $>1$, the moduli of abelian
varieties of CM-type and their torsion points generate abelian extensions, not
of the field of complex multiplication, but of another field associated with
it --- the latter is now called the reflex field. In principle, the theory of
Shimura and Taniyama allows one to list the abelian varieties of CM-type whose
reflex field is contained in a given CM-field $E$, and to determine the
extensions of $E$ obtained from the moduli. However, the results in the
published literature are unsatisfactory --- for example, they don't give a
good description of the largest abelian extension of a field obtainable in
this fashion (see Shimura and Taniyama 1961, Chapter IV; Shimura 1962; Shimura
1998, Chapter IV). Thus, the next theorem is of considerable interest.

\begin{theorem}[Wei
1993, 1994] Let $E$ be a CM-field.  Let $F$ be the maximal totally real subfield
of $E$, and let $H$ be the image of $\Gal(F^{\text{\rm ab}}/F\cdot
\mathbb{Q}^{\text{\rm ab}})$ in $\Gal(E^{\text{\rm ab}}/E)$ under the
Verlagerung map
\[\Gal(\mathbb{Q}^{\al}/F)^{\text{\rm ab}}\to\Gal(\mathbb{Q}^{\al}/E)^{\text{\rm ab}}
.\]
Then the field obtained by adjoining to $E$ the moduli of all polarized abelian
varieties of CM-type (and their torsion points) with reflex field contained
in $E$ is
\[{\cal M}_E=E^{\text {\rm ab}}{}^H\]
\end{theorem}

The theorem is proved by combining the three lemmas below.

Let $T$ be a torus over $\mathbb{Q}$ and $\mu$ a cocharacter of $T$. We are
only interested in pairs $(T,\mu)$ satisfying the conditions:

\begin{enumerate}
\item \ $T$ is split by a CM-field; equivalently, for all automorphisms $
\tau$ of $\mathbb{C}$, the actions of $\tau\iota$ and $\iota\tau$ on
$X^{*}(T)$ agree;

\item \ the weight $-\mu-\iota\mu$ of $\mu$ is defined over $\mathbb{Q}$.
\end{enumerate}

Let $(T,\mu)$ be a pair satisfying (a) and (b). Its \emph{reflex field} $
E(T,\mu)$ is the field of definition of $\mu$ --- because of (a), $E(T,\mu)$
is a subfield of a CM-field. Let $E\supset E(T,\mu)$. On applying
$\Res_{E/\mathbb{Q}}$ (Weil restriction) to the homomorphism $\mu
:\mathbb{G}_{m}\to T_{E}$ and composing with the norm map, we obtain a
homomorphism $N(T,\mu)$:
\[
\begin{CD}
\Res_{E/{\mathbb{Q}}}{\mathbb{G}}_m
@>{\Res_{E/{\mathbb{Q}}}\mu}>>
\Res_{E/{\mathbb{Q}}}T_E
@>\text{Norm}_{E/\mathbb{Q}}>>T
\end{CD}
\]
For any $\mathbb{Q}$-algebra $R,$ this gives a homomorphism
\[
(E\otimes_{\mathbb{Q}}R)^{\times}\to T(R).
\]
Let $\overline{T(\mathbb{Q})}$ be the closure of $T(\mathbb{Q})$ in $
T(\mathbb{A}_{f})$. The \emph{reciprocity map }
\[
r(T,\mu):\Gal(E^{\text{\textrm{\emph{ab}}}}/E)\to T(\mathbb{A}_{f})/\overline{
T(\mathbb{Q})}%
\]
is defined as follows: let $\tau\in\Gal(E^{\text{\textrm{ab}}}/E)$, and let $
t\in\mathbb{A}_{E}^{\times}$ be such that $\rec_{E}(t)=\tau$; write
$t=t_{\infty}\cdot t_{f}$ with $t_{\infty}\in(E\otimes_{\mathbb{Q}}%
\mathbb{R})^{\times}$ and $t_{f}\in(E\otimes_{\mathbb{Q}} \mathbb{A}%
_{f})^{\times}$; then
\[
r(T,\mu)(\tau)\overset{\text{\rm df}}{=}N(T,\mu)(t_{f})\mod\overline {
T(\mathbb{Q})}.
\]

\begin{lemma}
Let $E$ be a CM-field, and let $H$ be as in the statement of the theorem.
Then
\[H=\bigcap\Ker(r(T,\mu ))\]
where $(T,\mu )$ runs over the pairs satisfying (a) and (b) and such that
$E(T,\mu )\subset E$.
\end{lemma}
\begin{proof}
There is a universal such pair, namely, $(S^E,\mu^E)$, and so
\[\bigcap\Ker r(T,\mu )=\Ker r(S^E,\mu^E).\]
Because $S^E$ has no $\mathbb{R}$-split subtorus that is not already split over $
\mathbb{Q}$, $S^E(\mathbb{Q})$
is closed in $S^E(\mathbb{A}_f)$.  Thus, to prove the lemma, one must show that $
H$ is
the kernel of
\[r(S^E,\mu^E):\Gal(E^{\text {\rm ab}}/E)\to S^E(\mathbb{A}_f)/S^E(
\mathbb{Q}).\]
This can be done by direct calculation (Wei 1994, Theorem 2.1).
\end{proof}

For any CM-field $K$ with CM-type $\Phi$, we obtain a pair $(K^{\times}
,\mu_{\Phi})$ satisfying (a) and (b) (see \S 4).

\begin{lemma}
Let $E$ be a CM-field, and let $H$ be as above.  Then
\[H=\bigcap\Ker r(K^{\times},\mu_{\Phi})\]
where the intersection is over all CM-types $(K,\Phi )$ with reflex field
contained in $E$.
\end{lemma}
\begin{proof}
For each $(K,\Phi )$ with reflex field contained in $E$, we obtain a homomorphism
$\rho_{\Phi}:S^E\to K^{\times}$ (see \S 4), and (cf. the preceding proof) it suffices to show that
$\bigcap\Ker\rho_{\Phi}=1$.  But $X^{*}(S^E)$ is generated by the CM-types $
\Psi$ on $E$, and $\Psi$
occurs in $\rho_{\Phi}$ for $\Phi$ the reflex of $\Psi$ (ibid. 1.5.1).
\end{proof}

\begin{lemma}
Let $(A,i)$ be an abelian variety over $\mathbb{C}$ of CM-type $(K,
\Phi )$, and let $E$ be the
reflex field of $(K,\Phi )$.  The field of moduli of $(A,i)$ and its torsion points is
$(E^{\text {\rm ab}})^{H(\Phi )}$ where $H(\Phi )$ is the kernel of $
r(K^{\times},\Phi )$.
\end{lemma}
\begin{proof}
This is (yet another) restatement of the Theorem of Shimura and Taniyama.
\end{proof}

In fact, (ibid.) for a CM-field $E$, the following fields are equal:

\begin{enumerate}
\item \ the fixed field of $H$;

\item \ the field generated over $E$ by the fields of moduli of all CM-motives
and their torsion points with reflex field contained in $E$;

\item \ the field generated over $E$ by the fields of moduli of the CM-motive
and its torsion points defined by any faithful representation of $ S^{E}$;

\item \ the field generated over $E$ by the fields of moduli of the polarized
abelian varieties and their torsion points of CM-type with reflex field
contained in $E$;
\end{enumerate}

Moreover, for some Siegel modular variety and special point $z$, this is the
field generated by the values at $z$ of the $E$-rational modular functions on
the variety (ib. 3.3.2; see also the next section).

\subsection{Special Values of Modular Functions}

\footnote{This subsection is a manuscript of mine dated May 6, 1993.}

Abelian class field theory classifies the abelian extensions of a number field
$k$, but does not explain how to generate the fields. In his Jugendtraum,
Kronecker suggested that the abelian extensions of $\mathbb{Q}$ can be
generated by special values of the exponential function, and that the abelian
extensions of an imaginary quadratic number field can be generated by special
values of elliptic modular functions. This idea of generating abelian
extensions using special values of holomorphic functions was taken up by
Hilbert in his twelfth problem, where he suggested ``finding and discussing
those functions that play the part for any algebraic number field
corresponding to that of the exponential function for the field of rational
numbers and of the elliptic modular functions for imaginary quadratic number fields.''

Here we explain how the theory of Shimura varieties allows one to define a
class of modular functions naturally generalizing that of the elliptic modular
functions, and that it allows one to identify the fields generated by the
special values of the functions as the fields of moduli of CM-motives.

\subsubsection{Modular functions over $\mathbb{C}$.}

To define a Shimura variety, one needs a reductive group $G$ over $\mathbb{Q}$
and a $G(\mathbb{R)}$-conjugacy class $X$ of homomorphism $ \mathbb{S\to
G_{R}}$ satisfying the following conditions:

\begin{description}
\item[SV1] for each $h\in X$, the Hodge structure on the Lie algebra
$\mathfrak{g}$ of $G$ defined by $\Ad\circ h:\mathbb{S\to\GL(\mathfrak{g_{R}
)}}$ is of type $\{(-1,1),(0,0),(-1,1)\}$;

\item[SV2] for each $h\in X$, $\ad h(i)$ is a Cartan involution on $
G_{\mathbb{R}}^{\text{ad}}$;

\item[SV3] the adjoint group $G^{\text{ad}}$ of $G$ has no factor defined over
$\mathbb{Q}$ whose real points form a compact group, and the identity
component of the centre of $G$ splits over a CM-field.
\end{description}

The condition (SV1) implies that the restriction of $h$ to $\mathbb{G_{m}%
\subset S}$ is independent of $h\in X$. We denote its reciprocal by
$w_{X}:\mathbb{G_{m}\to G_{C}}$, and call it the \textit{weight\/} of the
Shimura variety. The weight is always defined over a totally real number
field, and we shall be especially interested in Shimura varieties for which it
is defined over $\mathbb{Q}$.

Consider a pair $(G,X)$ satisfying the Axioms (SV1-3). The set $ X$ has a
canonical $G(\mathbb{R)}$-invariant complex structure for which the connected
components are isomorphic to bounded symmetric domains.

For each compact open subset $K$ of $G(\mathbb{A}_{f})$,
\[
\text{Sh}_{K}(G,X)\overset{\text{\rm df}}{=}G(\mathbb{Q})\backslash X\times
G(\mathbb{A}_{f})/K
\]
is a finite disjoint union of quotients of the connected components of $ X$ by
arithmetic subgroups of $G^{\text{ad}}(\mathbb{Q})^{+}$, say
\[
\text{Sh}_{K}(G,X)=\bigcup\Gamma_{i}\backslash X_{i}.
\]
For $K$ sufficiently small, each space $\Gamma_{i}\backslash X_{i}$ will be a
complex manifold, and, according to Baily and Borel (1966), it has a natural
structure of a quasi-projective variety over $\mathbb{C}$. Hence
$\text{Sh}_{K}(G,X)$ is an algebraic variety over $\mathbb{C}$, and the
\textit{Shimura variety} $\text{Sh} (G,X)$ is the projective system of these
varieties, or (what amounts to the same thing) the limit of the system,
together with the action of $G(\mathbb{A_{f})}$ defined by the rule:
\[
[x,a]\cdot g=[x,ag],\quad x\in X,\quad a,g\in G(\mathbb{A_{f}). }%
\]

A rational function $f$ on $\text{Sh}_{K}(G,X)$ is called an
\textit{automorphic function over} $\mathbb{C}$ when $\dim X>0$\textit{.\/}
Such a function defines (for each $i$) a meromorphic function $f_{i}$ on each
$X_{i}$ invariant under $ \Gamma_{i}$. Conversely a family $(f_{i})$ of
invariant meromorphic functions defines an automorphic function $f$ provided
each $f_{i}$ is ``meromorphic at infinity'' (this condition is automatically
satisfied except when $X_{i}$ has dimension $1$).

When the weight $w_{X}$ of the Shimura variety is defined over $\mathbb{Q}$,
we shall call the automorphic functions \textit{modular functions.\/}
Classically this name is reserved for functions on Shimura varieties that are
moduli variety for abelian varieties, but it is known that most Shimura
varieties with rational weight are moduli varieties for abelian motives, and
it is hoped that they are all moduli varieties for motives, and so our
nomenclature is reasonable. This class of functions is the most natural
generalization of the class of elliptic modular functions.

Note that it doesn't yet make sense to speak of the algebraic (much less
arithmetic) properties of the special values of modular functions, because,
for example, the product of a modular function with a complex number is again
a modular function.

\begin{example}
Let $G=\GL_2$ and let $X$ be the
$G(\mathbb{R})$-conjugacy class of homomorphism $\mathbb{S}\to\GL_{2,\mathbb{R}}$ containing the
homomorphism
$$a+ib\mapsto\left(\begin{matrix} a&-b\\
b&a\end{matrix} \right).$$
The map $h\mapsto h(i)\cdot i$ identifies $X$ with $\{z\in \mathbb{C}
\mid\Re (z)\neq 0\}$, and in
this case $\text{Sh}_K(G,X)$ is a finite union of elliptic modular curves
over $\mathbb{C}$.  If $K=\GL_2(\hat{\mathbb{Z}})$, then the field of modular functions on
$\text{Sh}_K(G,X)$ is $\mathbb{C}[j]$.
\end{example}

\begin{example}
Let $T$ be a torus over $\mathbb{Q}$ split by a
CM-field, and let $\mu\in X_{*}(T)$.  Define $h:{\mathbb{S}}\to T_\mathbb{R}$ by $
h(z)=\mu (z)\cdot\overline {\mu (z)}$.
Then $(T,\{h\})$ defines a Shimura variety.
\end{example}

\begin{remark}
Shimura varieties have been studied for 200 years\footnote{This seems to be
an exaggeration.}...Gauss, Picard, Poincar\'e, Hilbert, Siegel, Shimura,...  The
axiomatic definition given above is due to Deligne (except that he doesn't
require that the identity component of the centre split over a CM-field).
The name is due\footnote{Earlier Shimura curves had been  so
named by Ihara.} to Langlands.
\end{remark}

\subsubsection{Special points}

A point $x\in X$ is said to be \textit{special\/} if there exists a torus
$T\subset G$ (this means $ T$ is rational over $\mathbb{Q}$), such that
$\im(h_{x})\subset T_{\mathbb{R}}$. By a \textit{special pair} $ (T,x)$ in
$(G,X)$ we mean a torus $T\subset G$ together with a point $x\in X$ such that
$h_{x}$ factors through $T_{\mathbb{R}}$.

\begin{example}In the Example 2, the special
points correspond to points $z\in \mathbb{C}\setminus \mathbb{R}$ such that $
[\mathbb{Q}[z]:\mathbb{Q}]=2$.  For
such a $z$, the choice of a $\mathbb{Q}$-basis for $E=_{df}\mathbb{Q}[z]$ determines an embedding
$\mathbb{Q}[z]^{\times}\hookrightarrow\GL_2(\mathbb{Q})$, and hence an embedding $
T=_{df}(\mathbb{G}_m)_{E/\mathbb{Q}}\hookrightarrow\GL_2$.  The map $h_z$
factors through $T_\mathbb{R}\hookrightarrow\GL_{2,\mathbb{R}}$.
\end{example}

\subsubsection{Modular functions defined over number fields}

To a torus $ T$ defined over $\mathbb{Q}$ and a cocharacter $\mu$ of $T$
defined over a number field $ E$, we attach a \emph{reciprocity map }
\[
\begin{CD}
r(T,\mu):\Gal(E^{\text{ab}}/E)@>>>T(\mathbb{A}_f)/T(\mathbb{Q})^{-}
\end{CD}
\]
as in (Milne 1992\emph{,\/} p164)\footnote{Better, see above}. The
\emph{reflex field } $E(G,X)$ is defined to be the field of definition of the
$G(\mathbb{C} )$-conjugacy class of homomorphisms $\mathbb{G}_{m}\to
\mathbb{G}_{\mathbb{C}}$ containing $\mu_{x}$ for $x\in X$. It is a number
field, and is a subfield of a CM-field. Hence it is either itself a CM-field
or is totally real.

By a model of $\text{Sh}(G,X)$ over a subfield $k$ of $\mathbb{C}$, we mean a
scheme $ S$ over $k$ endowed with an action of $G(\mathbb{A_{f})}$ (defined
over $k$) and a $G(\mathbb{A_{f})}$-equivariant isomorphism $\text{Sh}(G,X)\to
S\otimes_{k}\mathbb{C}$. We use this isomorphism to identify $\text{Sh}%
(G,X)(\mathbb{C)}$ with $S(\mathbb{C)}$.

\begin{theorem}
There exists a model of $\text{Sh}(G,X)$ over $
E(G,X)$
with the following property: for all special pairs $(T,x)\subset
(G,X)$ and
elements $a\in G(\mathbb{A}_f)$, the point $[x,a]$ is rational over $
E(T,x)^{\text{ab}}$ and
$\tau\in\Gal(E(T,x)^{\text{ab}}/E(T,x))$ acts on $[x,a]$ according to the rule:
$$\tau [x,a]=[x,ar(\tau )]\text{, where }r=r(T,\mu_x).$$
The model is uniquely determined by this condition up to a unique
isomorphism.
\end{theorem}

The model in the theorem is said to be \textit{canonical}.

\begin{remark}
For Shimura varieties of PEL-type,
models over number fields were constructed by Mumford and
Shimura (and his students Miyake and Shih).  That they satisfy
condition in the theorem follows from the theorem of Shimura and
Taniyama.  Shimura defined the notion of a canonical model more
generally, and proved the existence in one interesting case where
the weight is not defined over $\mathbb{Q}$.  Deligne modified Shimura's
definition, proved that the canonical model is unique (if it exists)
(1971), and showed that it exists for all Shimura varieties of abelian
type (1979a).  In 1981 Borovoi suggested using a trick of
Piateski-Shapiro to extend the proof to the remaining cases, and this
was carried out by Milne in 1982 (Milne 1983).
\end{remark}

Write $\text{Sh}(G,X)_{E}$ for the model in the theorem, and for any $
k\supset E$, write $\text{Sh}(G,X)_{k}$ for $\text{Sh}(G,X)_{E}\otimes_{E}k$.

For a connected variety $V$ over a field $k$, the field of rational functions
on $V$ is a subfield of the field of rational functions on $V\otimes
_{k}\mathbb{C}$. We say that a modular function $f$ on $\text{Sh}_{K}(G,X)$ is
\textit{rational over a subfield} $k$ of $\mathbb{C}$ if it arises from a
rational function on $\text{Sh}_{K}(G,X)_{k}$.

Let $x\in X$ be special, say $\im(h_{x})\subset T_{\mathbb{R}}$. Then the
field of definition of $\mu_{x}$ is written $E(x)$ --- it is the reflex field
of $(T,h_{x})$, and is a finite extension of $E(G,X)$.

\subsubsection{The fields generated by special values of modular functions}

Let $V$ be a connected algebraic variety over a field $k$. A point $P\in
V(k^{\text{al}})$ is a morphism $\Spec k^{\text{al}}\to V$ --- we also use $
P$ to denote the image of the map. It corresponds to a $k$-homomorphism
$\mathcal{O}_{ V,P}\to k^{\text{al}}$. This homomorphism factors through
$\mathcal{O}_{V,P}/\mathfrak{m_{P}}$, and hence its image in $k^{\text{al}}$
is a subfield of $k^{\text{al}}$, which we denote $ k[P]$.

For any open affine neighbourhood $U$ of $P$, the field of rational functions
$k(V)$ on $V$ is the field of fractions of $k[U]$. For $f=g/h\in k(V)$, we can
speak of
\[
f(P)\overset{\text{\rm df}}{=}g(P)/h(P)\in k^{\text{al}}%
\]
whenever $f$ does not have a pole at $P$, i.e., when $h\notin\mathfrak{m}_{
\mathfrak{P}}$.

\begin{lemma}With the above notations,
$$k[P]=\bigcup k[f(P)]$$
where the union runs over the $f\in k(V)$ without a pole at $P$ (i.e.,
over $f\in {\cal O}_{V,P}$).
\end{lemma}
\begin{proof}
We may replace $V$ with an open affine
neighbourhood, and embed $V$ in $\mathbb{A}^n$.  Then $k[P]$ is the field
generated by the coordinates $(a_1,\ldots ,a_n)$ of $P$.  Clearly, for any
rational function $f(X_1,\cdots ,X_n)$ with coordinates in $k$,
$f(a_1,\ldots ,a_n)\in k[P]$ (if it is defined).  Conversely,
$k[P]=\bigcup k[f(a_1,\ldots ,a_n)]$ where $f$ runs through the polynomials in
$X_1,\ldots ,X_n$.
\end{proof}

For a number field $k$, let $k_{c}$ be the subfield of $k^{\text{al}}$
corresponding to
\[
\bigcap\Ker(r(T,\mu))
\]
where $(T,\mu)$ runs over the pairs $(T,\mu)$ consisting of a torus $ T$ split
by a CM-field and $\mu$ is a cocharacter of $T$ whose weight $-\mu-\iota\mu$
is defined over $\mathbb{Q}$. (Equivalently over the pairs $(T,\mu)$
consisting of a torus $T$ over $ \mathbb{Q}$ and a cocharacter $\mu$ of $T$
satisfying the Serre condition\footnote{This is the condition $(\sigma
-1)(\iota+1)\mu=0=(\iota+1)(\sigma-1)\mu$.}.)

\begin{theorem}
Let $k$ be an algebraic number field.
For any Shimura variety $\text{Sh}(G,X)$ such that $E(G,X)\subset
k$,
modular function $f$ on $\text{Sh}_K(G,X)$ rational over $k$, and special point $
x$
of $X$ such that $E(x)\subset k$, $f(x)\in k_c$ (if it is defined, i.e.,  $
f$ doesn't
have a pole at x).  Moreover, if $k$ contains a CM-field, then $k_
c$ is
generated by these special values.
\end{theorem}
\begin{proof}
Let $(T,x)$ be a special pair in $(G,X)$.  I claim that $
T$
splits over a CM-field.  To prove this, it suffices to show that
the action $\iota$ on $X^{*}(T)$ (or even $X^{*}(T)\otimes \mathbb{Q}
)$ commutes with that of $\tau$,
for all $\tau\in\Gal(\mathbb{Q}^{\text{al}}/\mathbb{Q})$.  But
$$X^{*}(T)\otimes \mathbb{Q}=X^{*}(T')\otimes \mathbb{Q}\oplus X^{*}(G^{\text{ab}}
)$$
where $T'=T/Z(G)$ (use that $G\to G^{\text{ad}}\times G^{\text{ab}}$ is an isogeny).  By
assumption (SV3), $X^{*}(G^{\text{ab}})$ splits over a CM-field, and it follows
from (SV2) that $\iota$ acts as -1 on $X^{*}(T')$ and hence commutes with
everything.  From Theorem 5, it is clear that $k([x,1])$ is fixed by
$\Ker(r(T,\mu_x))$ and so is contained in $k_c$.  From the lemma, this
implies that $f(x)\in k_c$ for all $f$.
Before proving the converse, we need a construction.  Let $E$ be
a CM-field, with maximal totally real subfield $F$.  Let $N$ be the
kernel of
\[\begin{CD}
(\mathbb{G}_m)_{E/\mathbb{Q}}@>>>S^E.
\end{CD}\]
It is a subgroup of $(\mathbb{G}_m)_{F/\mathbb{Q}}$, and hence is contained in the centre
of $\GL_{2,F}$, and we define $G=\GL_{2,F}/N$.  The choice of a basis for
$E$ as an $F$-space determines an inclusion $(\mathbb{G}_m)_{E/\mathbb{Q}}
\hookrightarrow\GL_{2,F}$, and
hence an inclusion $S^E\hookrightarrow G$.  Let $X$ be the $G(\mathbb{R}
)$ conjugacy class of
the composite
\[\begin{CD}
\mathbb{S}@>{h_{\text{can}}}>>S^E@>>>G.
\end{CD}\]
Then $\text{Sh}(G,X)$ is a Shimura variety of dimension $[F:\mathbb{Q}
]$ with weight
defined over $\mathbb{Q}$ and whose reflex field is $\mathbb{Q}$.
On applying this construction to the largest CM-field contained in
$k$, we obtain a Shimura variety $\text{Sh}(G,X)$ containing $\text{Sh}
(S^k,h_{\text{can}})$,
where $S^k=S^E$ is the Serre group for $k$.  The statement is now
(more-or-less) obvious.
\end{proof}

\subsection{Nonabelian solutions to Hilbert's Twelfth Problem}

By applying the new Main Theorem of Complex Multiplication (Theorem 9.1) in
place of the original, one obtains explicit non-abelian extensions of number
fields (Milne and Shih 1981, \S 5).

\section{Algebraic Hecke Characters are Motivic}

\footnote{This section is the notes of a seminar talk.} Weil's Hecke
characters of type $A_{0}$ are now called algebraic Hecke characters. In this
section, I show that they are all motivic (and explain what this means).

\subsection{Algebraic Hecke characters}

In this subsection, I explain the description of algebraic Hecke characters
given in Serre 1968.

\emph{Notations:} $\mathbb{Q}^{\al}$ is the algebraic closure of $\mathbb{Q}$
in $ \mathbb{C}$; $K$ is a fixed CM-field, $\Sigma=\Hom(K,\mathbb{Q}%
^{\al})=\Hom(K,\mathbb{C})$, and $\mathbb{I}= \mathbb{I}_{\infty}%
\times\mathbb{I}_{f}$ is the group of id\`eles of $ K$. For a finite extension
$k^{\prime}/k$, $(\mathbb{G}_{m})_{k^{\prime}/k}$ is the torus over $ k$
obtained from $\mathbb{G}_{m/k^{\prime}}$ by restriction of scalars.

I define an \emph{algebraic Hecke character\/} to be a continuous homomorphism
$\chi:\mathbb{I}\to\mathbb{Q}^{\al}{}^{\times}$ such that

\begin{enumerate}
\item $\chi=1$ on $\mathbb{I}_{\infty}$;

\item the restriction of $\chi$ to $K^{\times}\subset\mathbb{I}$ is given by
an algebraic character of the torus $(\mathbb{G}_{m})_{K/\mathbb{Q}}$.
\end{enumerate}

Condition (b) means that there exists a family of integers $(n_{\sigma}
)_{\sigma\in\Sigma}$ such that $\chi(x)=\prod\sigma(x)^{n_{\sigma}}$ for all
$x\in K^{\times}\subset\mathbb{I}$. Condition (a) means that $\chi$ factors
through $\mathbb{I}\to\mathbb{I}_{f}$. Thus, there is a one-to-one
correspondence between algebraic Hecke characters and continuous homomorphisms
$\mathbb{I}_{f}\to\mathbb{Q}^{\al\times}$ satisfying the analogue of (b)
(restriction to $K^{\times}\subset\mathbb{I}_{f}$) (cf. the definition in
Harder and Schappacher).

Let $\chi$ be a Hecke character.

\subsubsection{The character $\chi$ admits a modulus}

Let $\mathfrak{m}$ be a modulus for $ K$. Because $K$ has no real primes,
$\mathfrak{m}$ can be regarded as an integral ideal $\prod_{v}\mathfrak{p}%
_{v}^{m_{v}}$. Define
\[
W_{\mathfrak{m}}=\prod_{v|\infty}K_{v}^{\times}\times\prod_{v|\mathfrak{m}}
(1+\hat{\mathfrak{p}}_{v}^{m_{v}})\times\prod U_{v}%
\]
(as in my class field theory notes, Milne 1997, about V.4.6). The $
W_{\mathfrak{m}}$'s are open subgroups of $\mathbb{I}$ and any neighbourhood
of $1$ containing $ \mathbb{I}_{\infty}$ contains a $W_{\mathfrak{m}}$. Let
$V$ be a neighbourhood of $1$ in $\mathbb{C}^{\times}$ not containing any
subgroup $ \neq1$. Because $\chi$ is continuous and $1$ on $\mathbb{I}%
_{\infty}$, $\chi(W_{\mathfrak{m}})\subset V$ for some $\mathfrak{m}$, and
hence $\chi(W_{\mathfrak{m}})=1$. Such an $\mathfrak{m}$ will be called a
\emph{modulus\/} for $ \chi$. If $\mathfrak{m}$ is a modulus for $\chi$ and
$\mathfrak{m}|\mathfrak{m}^{\prime}$, then $\mathfrak{m}^{\prime}$ is also a
modulus for $ \chi$.

\subsubsection{The infinity type of $\chi$}

Let $\mathbb{Z}^{\Sigma}$ be the free abelian group generated by $\Sigma$,
with $\tau\in\Gal(\mathbb{Q}^{\al}/\mathbb{Q})$ acting by $ \tau(\sum
_{\sigma\in\Sigma}n_{\sigma}\sigma)=\sum n_{\sigma}(\tau\circ\sigma)$. The
character group $X^{*}((\mathbb{G}_{m})_{K/\mathbb{Q}})=\mathbb{Z}^{\Sigma}$,
and so $ \chi|K^{\times}=\sum n_{\sigma}\sigma$ for some $n_{\sigma}%
\in\mathbb{Z}$. The element $\sum_{\sigma\in\Sigma}n_{ \sigma}\sigma$ is
called the \emph{infinity type\/} of $\chi$.

Let $U_{\mathfrak{m},1}=K^{\times}\cap W_{\mathfrak{m}}$ --- this is a
subgroup of finite index in the units $ U$ of $K$ defined by congruence
conditions at the primes dividing $\mathfrak{m}$. If $ \mathfrak{m}$ is a
modulus for $\chi$, $\chi=1$ on $U_{\mathfrak{m},1}$, and this implies that $
n_{\sigma}+n_{\bar{\sigma}}=\text{\textrm{constant}}$, independent of $\sigma$
(apply the Dirichlet unit theorem).

\subsubsection{The Serre group.}

Let $\Xi$ be the group of infinity types, i.e.,
\[
\Xi=\{\sum_{\sigma\in\Sigma}n_{\sigma}\sigma\mid n_{\sigma}+n_{ \bar{\sigma}%
}=\text{\textrm{constant}}\}\subset\mathbb{Z}^{\Sigma}.
\]
It is a free $\mathbb{Z}$-module of finite rank on which $\Gal(\mathbb{Q}%
^{\al} /\mathbb{Q})$ acts, and we define the \emph{Serre group} $S^{K}$ to be
the torus over $\mathbb{Q}$ with character group $ \Xi$. Thus, for any field
$L\subset\mathbb{Q}^{\al}$,
\[
S^{K}(L)=\Hom(X^{*}(S^{K}),\mathbb{Q}^{\al})^{\Gal(\mathbb{Q}^{\al}/L)} .
\]
Because $X^{*}(S^{K})\subset X^{*}((\mathbb{G}_{m})_{K/\mathbb{Q}})$, $S^{K}$
is a quotient of $(\mathbb{G}_{m})_{K/\mathbb{Q}}$. The map on $\mathbb{Q}%
$-rational points $K^{\times}\to S^{K}(\mathbb{Q})$ sends $x \in K^{\times}$
to the map $\xi\mapsto\xi(x)$, $\xi\in\Xi$.

\subsubsection{Serre's extension}

I claim that there exists a modulus $ \mathfrak{m}$ such that $U_{\mathfrak{m}%
,1}$ is contained in $\Ker(\xi)$ for all $\xi\in\Xi$. Indeed, in order for
$\xi=\sum n_{\sigma}\sigma$ to lie in $\Xi$, its restriction to the totally
real subfield $F$ of $ K$ must be a power of the norm. Thus all $\xi=1$ on
some subgroup $U$ of index at most $ 2$ in $U_{F}$. But $U_{F}$ is of finite
index in $U_{K}$ (Dirichlet unit theorem again), and so $U$ has finite index
in $U_{K}$. An old theorem of Chevalley states that every subgroup of finite
index in $U_{K}$ is a congruence subgroup, i.e., contains $U_{\mathfrak{m},1}$
for some $\mathfrak{m}$.

\emph{From now on,} $\mathfrak{m}$ \emph{will denote a modulus with this
property; thus the } \emph{canonical map} $K^{\times}\to S^{K}(\mathbb{Q})$
\emph{factors through} $ K^{\times}/U_{\mathfrak{m},1}$\emph{.}

Recall (e.g., Milne 1997, V.4.6) that $\mathbb{I}/W_{\mathfrak{m}}\cdot
K^{\times}=C_{\mathfrak{m}}$, the ray class group with modulus $\mathfrak{m}$
($=I^{S(\mathfrak{m})}/i(K_{\mathfrak{m},1})$). In particular, it is finite.
There is an exact sequence
\[
1\to K^{\times}/U_{\mathfrak{m},1}\to\mathbb{I}/W_{\mathfrak{m}}\to C_{
\mathfrak{m}}\to1.
\]

Serre shows that there is a canonical exact sequence of commutative algebraic
groups over $\mathbb{Q}$
\[
1\to S^{K}\to T_{\mathfrak{m}}\to C_{\mathfrak{m}}\to1
\]
(here $C_{\mathfrak{m}}$ is regarded as a finite constant algebraic group) for
which there is a commutative diagram
\[%
\begin{array}
[c]{cccccccccc}%
1 & \to & K^{\times}/U_{\mathfrak{m},1} & \to & \mathbb{I}/W_{\mathfrak{m}} &
\to & C_{ \mathfrak{m}} & \to & 1 & \\
&  & \downarrow &  & \downarrow\varepsilon &  & \| &  &  & \\
1 & \to & S^{K}(\mathbb{Q}) & \to & T_{\mathfrak{m}}(\mathbb{Q}) & \to &
C_{\mathfrak{m}} & \to & 1. &
\end{array}
\]
Moreover, there is a natural one-to-one correspondence between the algebraic
Hecke characters $\chi$ of $K$ admitting $\mathfrak{m}$ as a modulus and the
characters of $T_{\mathfrak{m}}$ as an algebraic group. [The proofs of these
statements are straightforward.] The algebraic Hecke character corresponding
to a character $\chi$ of $T_{\mathfrak{m}}$ is the composite
\[
\mathbb{I}\to\mathbb{I}/W_{\mathfrak{m}}\xr{\varepsilon}T_{\mathfrak{m}}
(\mathbb{Q}^{\al})\to\mathbb{Q}^{\al\times}.
\]
>From now on, I'll \emph{define\/} an algebraic Hecke character to be a
character of $T_{\mathfrak{m}}$ for some $\mathfrak{m}$. Its infinity type is
its restriction to $ S^{K}$. The Dirichlet characters are the Hecke characters
with trivial infinity type (and hence factor through $C_{\mathfrak{m}}$).

\emph{Warning!\/} Our notations differ from those of Serre---in particular, he
switches the $S$ and the $T$.

\subsubsection{The $\ell$-adic representation}

One checks that the two maps
\[
\alpha_{\ell}:\mathbb{I}\xr{\text{\rm proj}}(K\otimes_{\mathbb{Q}}\mathbb{Q}_{
\ell})^{\times}\to S^{K}(\mathbb{Q}_{\ell})\to T_{\mathfrak{m}}(\mathbb{Q}_{
\ell})
\]
and
\[
\varepsilon:\mathbb{I}\to T_{\mathfrak{m}}(\mathbb{Q})
\]
coincide on $K^{\times}$. Therefore, $\varepsilon_{\ell}\overset
{\text{\textrm{df}}}{ =}\varepsilon\cdot\alpha_{\ell}^{-1}:\mathbb{I}\to
T_{\mathfrak{m}}(\mathbb{Q}_{ \ell})$ factors through $\mathbb{I}/K^{\times
}\mathbb{I}_{\infty}$, and hence through $\Gal(K^{\ab}/K)$ --- thus
$\varepsilon_{\ell}$ is a continuous homomorphism
\[
\Gal(K^{\text{\rm ab}}/K)\to T_{\mathfrak{m}}(\mathbb{Q}_{\ell}).
\]

\subsubsection{The Hecke character in the usual sense.}

The same argument with $\ell$ replaced by $\infty$ gives a homomorphism
$\varepsilon_{ \infty}:\mathbb{I}\to T_{\mathfrak{m}}(\mathbb{R})$ that is
\emph{not } usually trivial on the connected component of $\mathbb{I}$. Its
composite with any character of $T_{\mathfrak{m}}$ defined over $\mathbb{C}$
is a Hecke character in the usual (broad) sense: continuous homomorphism
$\mathbb{I}\to\mathbb{C}^{\times}$ trivial on $ K^{\times}$.

\subsection{Motivic Hecke Characters}

Let $k$ be a subfield of $\mathbb{C}$. An abelian variety $A$ over $k$ is said
to be of \emph{CM-type\/} if there exists a product of fields $E\subset
\End^{0}(A_{\mathbb{C}})$ such that $H^{1}_{B}(A,\mathbb{Q} )$ is a free
$E$-module of rank $1$. It is said to have \emph{CM over} $k$ if
$E\subset\End^{0}(A)$. It is possible to choose $ E$ so that it is stabilized
by the Rosati involution of some polarization of $ A$, which implies that it
is a product of CM-fields.

Let $A$ be an abelian variety with CM by $E$ over $K$. Then $\Gal (K^{\al}/K)$
acts on $V_{\ell}(A)$ by $E\otimes_{\mathbb{Q}}\mathbb{Q}_{\ell}$-linear maps.
But $ V_{\ell}(A)$ is a free $E\otimes_{\mathbb{Q}}\mathbb{Q}_{\ell}$-module
of rank one, so this action defines a homomorphism
\[
\rho_{\ell}:\Gal(K^{\ab}/K)\to(E\otimes_{\mathbb{Q}}\mathbb{Q}_{\ell}
)^{\times}\subset\GL(V_{\ell}(A)).
\]
The main theorem of Shimura-Taniyama theory can be stated\footnote{Cf. 13.1.}
as follows: \begin{E}
For ${\mathfrak{m}}$ sufficiently large, there exists a unique homomorphism
$\chi :T_{{\mathfrak{m}}}\to (\mathbb{G}_m)_{E/\mathbb{Q}}$ of tori such $
\rho_{\ell}=\chi (\mathbb{Q}_{\ell})\circ\varepsilon_{\ell}$.
\end{E}

An embedding $\sigma$ of $E$ into $\mathbb{Q}^{\al\times}$, defines a
character $ \chi_{\sigma}$ of $T_{\mathfrak{m}}$, which (by definition) is an
algebraic Hecke character. Such characters are certainly motivic.

The infinity type of a Hecke character arising in this way is a CM-type on
$K$, i.e., $n_{\sigma}\geq0$, $n_{\sigma}+n_{\bar{\sigma}}=1$, and Casselman
showed that conversely every Hecke character with infinity type a CM-type
arises in this fashion.

More generally, I discussed motives of type $M=(A,e)$, $A$ an abelian variety
of CM-type, $e^{2}=e$, $e\in C^{g}(A\times A)/\!\!\sim$ (algebraic classes of
codimension $g=\dim A$ modulo numerical equivalence). Such an $M$ has an
endomorphism ring and Betti and \'etale cohomology groups, and so one can make
the same definitions as for $A$. Note that $M=(A,e)$ may have CM over $k$
without $ A$ having CM over $k$. The analogue of (15.1) holds. A Hecke
character arising from such a motive\emph{,\/} or the product of such a
character with a Dirichlet character, will be called \emph{motivic. }

If we assume the Hodge conjecture, then every algebraic Hecke character is motivic.

After a theorem of Deligne (1982a), we no longer need to assume the Hodge
conjecture, but at the cost of replacing $e$ with an absolute Hodge class.

\subsection{The proof}

(that all algebraic Hecke characters are motivic). The CM-motives discussed
above over a field $k$ form a category $\CM( k)$ that looks like the category
of representations of an algebraic group: it is $\mathbb{Q}$-linear, abelian,
has a tensor product, duals, and every object has rank equal to a nonnegative
integer. The theory of Tannakian categories then shows that it \emph{is\/} the
category of representations of a pro-algebraic group. (We are using absolute
Hodge classes to define motives.)

What is the pro-algebraic group? When $k=\mathbb{C}$, one sees easily that it
$S=\plim S^{K}$ (projective limit over the CM-subfields of $\mathbb{Q}^{\al}%
$). Hint: The abelian varieties of CM-type over $\mathbb{C}$ are classified up
to isogeny by CM-types, and $S^{K}$ is generated by the CM-types on $K$.

When $k=\mathbb{Q}^{\al}$, the group is again $S$ (base change $\mathbb{Q}%
^{\al} \to\mathbb{C}$ gives an equivalence of categories of CM-motives).

When $k=\mathbb{Q}$, the general theory tells us it is an extension
\[
1\to S\to T\to\Gal(\mathbb{Q}^{\al}/\mathbb{Q})\to1.
\]
Following Langlands, we call $T$ the Taniyama group.

Deligne, Grothendieck, and Serre knew in the 1960s that the general theory
predicted the existence of such an extension, but couldn't guess what it was.
(Although he doesn't say so, these ideas must have suggested to Serre his
interpretation of algebraic Hecke characters.) In the late 1970s, in trying to
understand the zeta functions of Shimura varieties, Langlands wrote down some
cocycles, which Deligne recognized should give the above extension. He
verified they do by proving that there is only one such extension having
certain natural properties shared by both extensions.

When $k=K\subset\mathbb{Q}^{\al}$, the group attached to the category of
CM-motives over $K$ is the subextension
\[
1\to S\to T^{K}\to\Gal(\mathbb{Q}^{\al}/K)\to1
\]
of the above extension. Thus, to give a CM-motive over $K$ is to give a
representation of $T^{K}$ on a finite-dimensional $\mathbb{Q}$-vector space.

>From Langlands's description of this extension, one sees that, for any $
\mathfrak{m}$, there is a canonical map (see 10.2a) of extensions:
\[%
\begin{array}
[c]{ccccccccc}%
1 & \to & S & \to & T^{K} & \to & \Gal(\mathbb{Q}^{\al}/K) & \to & 1\\
&  & \downarrow &  & \downarrow &  & \downarrow &  & \\
1 & \to & S^{K} & \to & T_{\mathfrak{m}} & \to & C_{\mathfrak{m}} & \to & 1
\end{array}
\]
Let $E$ be a CM-field. A homomorphism $T_{\mathfrak{m}}\to(\mathbb{G}%
_{m})_{E/\mathbb{Q}}$ defines by composition a representation $T^{K}%
\to(\mathbb{G}_{m})_{E/\mathbb{Q}}\hookrightarrow\GL (E^{\prime})$ where
$E^{\prime}=E$ regarded as a $\mathbb{Q}$-vector space. Therefore a Hecke
character $\chi$ defines a CM-motive $M (\chi)$ with CM by $E$ over $K$. The
motive $M(\chi)$ is related to $\chi$ as in 15.1, and so $ \chi$ is motivic.

\section{Periods of Abelian Varieties of CM-type}

Deligne's theorem (Deligne 1978, Deligne 1982a) allows one to define a
category of CM-motives over any field of characteristic zero (Deligne and
Milne 1982, \S 6).

Let $M$ be a simple CM-motive over $\mathbb{Q}^{\al}\subset\mathbb{C}$. Then
$\End (M)$ is a CM-field $K$. The Betti realization $H_{B}(M)$ of $M$ is a
vector space of dimension $ 1$ over $K$, and the de Rham realization
$H_{\dR}(M)$ is free of rank $1$ over $ K\otimes_{\mathbb{Q}}\mathbb{Q}^{\al}%
$. For $\sigma:K\hookrightarrow\mathbb{Q}^{\al}$, let $H_{\dR}(M)_{\sigma}$
denote the $ \mathbb{Q}^{\al}$-subspace of $H_{\dR}(M)$ on which $x\in K$ acts
as $\sigma(x)\in\mathbb{Q}^{\al}$. Then $H_{\dR}(M)$ being free of rank $ 1$
means that each $H_{\dR}(M)_{\sigma}$ has dimension $1$ and
\[
H_{\dR}(M)=\oplus_{\sigma:K\hookrightarrow\mathbb{Q}^{\al}}H_{\dR}
(M)_{\sigma}.
\]
Let $e$ be a nonzero element of $H_{B}(M)$, and let $\omega_{\sigma}$ be a
nonzero element of $H_{\dR}(M)_{\sigma}$. Under the canonical isomorphism
\[
H_{B}(M)\otimes_{\mathbb{Q}}\mathbb{C}\to H_{\dR}(M)\otimes_{\mathbb{Q}^{\al}}
\mathbb{C}%
\]
$e$ maps to a family $(e_{\sigma})$, $e_{\sigma}\in H_{\dR}(M)_{\sigma}
\otimes_{\mathbb{Q}^{\al}}\mathbb{C}$. Define $p(M,\sigma)\in\mathbb{C}$ by
the formula
\[
p(M,\sigma)\cdot e_{\sigma}=\omega_{\sigma}.
\]
When regarded as an element of $\mathbb{C}^{\times}/\mathbb{Q}^{\times}$, $
p(M,\sigma)$ is independent of the choices of $e$ and of $\omega_{\sigma}$ ---
the $p(M,\sigma)$ are called the \emph{periods\/} of $ M$. Clearly,
$p(M,\sigma)$ depends only on the isomorphism class of $M$.

Let $S$ be the Serre group --- it is the projective limit of the Serre groups
$S^{E}$ for $E$ a CM-field contained in $\mathbb{C}$. It is the protorus over
$ \mathbb{Q}$ whose character group $X^{*}(S)$ consists of locally constant
functions $\phi:\Gal(\mathbb{Q}^{\cm}/\mathbb{Q})\to\mathbb{Z}$ such that
$\phi(\tau)+\phi(\iota\tau)$ is independent of $\tau$. The Betti fibre functor
defines an equivalence from the category of CM-motives over $\mathbb{Q}^{\al}$
to the category of finite-dimensional representations of $ S$. Thus the set of
simple isomorphism classes of CM-motives over $\mathbb{Q}^{\al}$ is in natural
one-to-one correspondence with the set of $\Gal(\mathbb{Q}^{\al}/\mathbb{Q}
)$-orbits in $X^{*}(S)$. Let $\phi\in X^{*}(S)$, and let $M(\phi)$ be the
CM-motive corresponding to $ \phi$. The endomorphism algebra of $M(\phi)$ is
$K=\mathbb{Q}^{\al H}$, where $ H$ is the stabilizer of $\phi$ in
$\Gal(\mathbb{Q}^{\al}/\mathbb{Q})$. Thus, for each $\phi\in X^{*} (S)$ and
coset representative for $H$ in $\Gal(\mathbb{Q}^{\al}/\mathbb{Q})$, we obtain
a period
\[
p(\phi,\sigma)\overset{\text{\textrm{df}}}{=}p(M(\phi),\sigma).
\]
Any relation among the $\phi$'s yields an isomorphism among the motives, and
hence a relation among the periods. When $\phi$ is taken to be a CM-type, then
$p(\phi,\sigma)$ is the period of an abelian variety of CM-type. Thus, we see
that Deligne's theorem (Deligne 1978) yields an array of relations among the
periods of abelian varieties of CM-type. (See Deligne's talk at the Colloq.,
\'Ecole Polytech., Palaiseau, 1979 (Deligne 1980); also Shimura's talk at the
same conference (Shimura 1980).)

In this context, one should also mention Blasius 1986.

\section{Review of: Lang, Complex Multiplication, Springer 1983.}

The \footnote{This is the author's version of MR 85f:11042.} theory of complex
multiplication for elliptic curves describes how an automorphism of
$\mathbb{C}$ acts on an elliptic curve with complex multiplication and its
torsion points. As a consequence, when the curve is defined over a number
field, one obtains an expression for its zeta function in terms of Hecke
$L$-series. The theory was generalized to abelian varieties in so far as it
concerned automorphisms fixing the reflex field by Shimura, Taniyama, and Weil
in the fifties. As a consequence, when the abelian variety is defined over a
number field containing the reflex field, they obtained an expression for its
zeta function (except for finitely many factors) in terms of Hecke $L$-series.
A thorough account of this is given in Shimura and Taniyama (1961).
Improvements are to be found in Shimura 1971 (Sections 5.5 and 7.8). Serre and
Tate (1968) extended the result on the zeta function to all the factors, and
computed the conductor of the variety. Serre (1968), Chapters 1 and 2,
re-interpreted some of this work in terms of algebraic tori. In 1977 Langlands
made a conjecture concerning Shimura varieties which was shown to have as a
corollary a description of how every automorphism of $\mathbb{C}$ acts on an
abelian variety with complex multiplication and its torsion points, and in
1981 Deligne proved the corollary (Deligne et al. 1982). Since this gives an
expression for the zeta function of such a variety over any number field in
terms of Weil $L$-series, it completes the generalization to abelian varieties
of the basic theory of complex multiplication for elliptic curves.

The first four chapters of the Lang's book are devoted to the same material as
that in (the sections of) the works of Shimura and Taniyama, Shimura, and
Serre and Tate cited above: the analytic theory of abelian varieties with
complex multiplication, the reduction of abelian varieties, the main theorem
of complex multiplication, and zeta functions. Lang's account is less detailed
but probably more readable than his sources. For example, whereas Shimura and
Taniyama's discussion of reduction is painfully detailed (they, like the
author, use the language of Weil's Foundations), that of the author is brief
and sketchy. The result of Serre and Tate on the conductor is not included
and, in the statement of the main theorem, it is unnecessarily assumed that
the abelian variety is defined over a number field.

Chapter 5 discusses fields of moduli and the possibility of descending abelian
varieties with complex multiplication to smaller fields (mainly work of
Shimura), and Chapter 6 introduces some of the algebraic tori associated with
abelian varieties having complex multiplication and uses them to obtain
estimates for the degrees of the fields generated by points of finite order on
the varieties.

The final chapter (based on a manuscript of Tate)\footnote{This is what the
book says, but in fact \S 4 of Chapter 7 is essentially a translation of
Deligne 1981. In the same chapter, Lang (p175) credits a theorem of Deligne
and Langlands to Deligne alone and (p163, p171) a theorem of mine to Deligne
(Langlands's conjecture; joint with Shih for Shimura varieties of abelian type
and with Borovoi in general).} gives the most down-to-earth statement of the
new main theorem of complex multiplication (the corollary of Langlands's
conjecture) and includes part of the proof (but, unfortunately, only the more
technical, less illuminating, part). Zeta functions are not discussed in this
general context.

The exposition is very clear in parts, but in others it is marred by
carelessness. For example, in Chapter 3, the definition of $\mathfrak{a}%
$-multiplication is incorrect (the universal property is not universal), in
the proof of (3.1) it is nowhere shown that the reduction of an $\mathfrak{a}%
$-multiplication is an $\mathfrak{a}$-multiplication, and in the proof of the
Main Theorem 6.1 it is not possible to write the id\`ele $s$ in the way the
author claims on p. 82 under his assumptions.

In summary, this book will be useful, in much the same way as a good lecture
course, for someone wishing to obtain a first understanding of the subject,
but for a more complete and reliable account it will be necessary to turn to
the original sources mentioned in this review. \hfill\emph{James Milne} (1-MI)

\centerline{\textbf{Additional References}}

Baily, W., and Borel, A., Compactification of arithmetic quotients of bounded
symmetric domains. Ann. of Math. (2) 84, 1966, pp 442--528.

Blasius, Don, On the critical values of Hecke $L$-series, Ann. of Math. (2)
\textbf{124} (1986), no.~1, 23--63.

Deligne, Pierre, Travaux de Shimura, in \textit{S\'eminaire Bourbaki, 23\`eme
ann\'ee (1970/71), Exp. No. 389}, 123--165. Lecture Notes in Math., 244,
Springer, Berlin, 1971.

Deligne, P., Cycles de Hodge sur les vari\'et\'es ab\'eliennes, 4pp, 8th
January, 1978.

Deligne, Pierre, Cycles de Hodge absolus et p\'eriodes des int\'egrales des
vari\'et\'es ab\'eliennes, M\'em. Soc. Math. France (N.S.) \textbf{1980},
no.~2, 23--33.

Deligne, P., Letter to Tate, 8th October 1981.

Deligne, P., and Milne, J., Tannakian categories. In: Deligne et al. 1982, pp 101--228.

Milne, J., The action of an automorphism of $\mathbf{C}$ on a Shimura variety
and its special points. Arithmetic and geometry, Vol. I, 239--265, Progr.
Math., 35, Birkh\"{a}user Boston, Boston, Mass., 1983.

Milne, J., The points on a Shimura variety modulo a prime of good reduction.
In: The zeta functions of Picard modular surfaces, 151--253, Univ.
Montr\'{e}al, Montreal, PQ, 1992.

Milne, J., Class Field Theory, 222 pages, 1997. Available at
\texttt{www.math.lsa.umich.edu/$\sim$jmilne/}.

Milne, J., and Shih, Kuang-yen, Automorphism groups of Shimura varieties and
reciprocity laws, Am. J. Math. 103 (1981), 911--935.

Serre, J.-P., and Tate, J., Good reduction of abelian varieties, Ann. of Math.
(2) 88 (1968), 492--517.

Shimura, Goro, On the class-fields obtained by complex multiplication of
abelian varieties, Osaka Math. J. \textbf{14} (1962), 33--44.

Shimura, G., On abelian varieties with complex multiplication, Proc. LMS 34
(1977), 65--86.

Shimura, G., The periods of abelian varieties with complex multiplication and
the special values of certain zeta functions, M\'em. Soc. Math. France (N.S.)
\textbf{1980}, no.~2, 103--106.

Shimura, G., Abelian Varieties with Complex Multiplication and Modular
Functions, Princeton, 1998.

Shimura, G., and Taniyama, Y., Complex multiplication of abelian varieties and
applications to number theory, Math. Soc. Japan, Tokyo, 1961.

Wei, W., Weil numbers and generating large field extensions. Thesis,
University of Michigan, 1993.

Wei, W., Moduli fields of CM-motives applied to Hilbert's 12-th problem, 17
pages, May 18, 1994.

Weil, A., On a certain type of characters of the id\`ele-class group of an
algebraic number-field, in \emph{Proceedings of the international symposium on
} \emph{algebraic number theory, Tokyo \& Nikko, 1955}, 1--7., Science Council
of Japan, Tokyo, 1956a.

Weil, A., On the theory of complex multiplication, in \emph{Proceedings of the
} \emph{international symposium on algebraic number theory,Tokyo \& Nikko,
1955}, 9--22., Science Council of Japan, Tokyo, 1956b.

Weil, A., Y. Taniyama (lettre d'Andr\'e Weil). Sugako-no Ayumi. 6 (1959), 21--22.

\end{document}

%% file: defs.tex
\DeclareMathSymbol{:}{\mathpunct}{operators}{"3A}

\numberwithin{equation}{section}

\def\1{{1\mkern-7mu1}}  
\newcommand\ad{\operatorname{ad}}
\newcommand\Ad{\operatorname{Ad}}
\newcommand\Aut{\operatorname{Aut}}

\newcommand\End{\operatorname{End}}

\newcommand\Gal{\operatorname{Gal}}
\newcommand\GL{\operatorname{GL}}

\newcommand\Hom{\operatorname{Hom}}

\newcommand\id{\operatorname{id}}

\newcommand\im{\operatorname{Im}}  

\newcommand\Ker{\operatorname{Ker}}

\newcommand\Nm{\operatorname{Nm}}

\let\plim=\varprojlim

\newcommand\rec{\operatorname{rec}}

\newcommand\Res{\operatorname{Res}}

\newcommand\Spec{\operatorname{Spec}}

\newcommand\Tr{\operatorname{Tr}}

\def\CM{{\bold{CM}}}

\def\ab{{\text{ab}}}
\def\al{{\text{al}}}
\def\ad{{\text{ad}}}

\def\cm{{\text{cm}}}

\def\dR{{\text{dR}}}

\let\xr=\xrightarrow